\documentclass[11pt]{article}
\usepackage[utf8]{inputenc}
\usepackage{graphicx}
\usepackage{caption}
\usepackage{amsmath}
\usepackage{amssymb}
\usepackage{float}
\usepackage{hyperref}
\usepackage{amsthm}
\usepackage{xcolor}
\usepackage{comment}
\usepackage{amsfonts}
\usepackage{dsfont}
\usepackage{stmaryrd}
\usepackage{xcolor}
\usepackage{subcaption}

\usepackage{comment} \usepackage{siunitx}
\usepackage{mathtools} 
\mathtoolsset{showonlyrefs} 

\usepackage{hyperref}

\DeclareMathOperator{\R}{\mathbb{R}}
\def\Ac{{\cal A}}
\def\Dc{{\cal D}}
\def\Lc{{\cal L}}
\def\Ec{{\cal E}}

\def\Kc{{\cal K}}

\def\Lc{{\cal L}}

\def\Nc{{\cal N}}
\def\Pc{{\cal P}}
\def\Sc{{\cal S}}
\def\Uc{{\cal U}}
\def\Tc{{\cal T}}
\def\Vc{{\cal V}}

\def\Yc{{\cal Y}}
\def\Wc{{\cal W}}
\def\Zc{{\cal Z}}
\def \E{\mathbb{E}}

\def \P{\mathbb{P}}

\def\1{{\bf 1}}

\def \N{\mathbb{N}}

\def \Sum{\displaystyle\sum}

\def \Frac{\displaystyle\frac}

\def \b1{\bf{1}}

\def\Dx#1{\Frac{\partial #1}{\partial x}}
\def\Dxx#1{\Frac{\partial^2 #1}{\partial x^2}}
\def\Dt#1{\Frac{\partial #1}{\partial t}}

\def\bop{{\boldsymbol p}}
\def\bopi{{\boldsymbol \pi}}

\def\botheta{{\boldsymbol \theta}}

\def\bolx{{\boldsymbol x}}

\def \mrb{\mathrm{b}} 
\def \d{\mathrm{d}} 
 
\def \mrp{\mathrm{p}} 
 
\def \mrt{\mathrm{t}}

\def\blue[#1]{{\color{blue}#1}}
\def\rd#1{#1}
\def\red#1{#1}

\def\beqs{\begin{eqnarray*}}
\def\enqs{\end{eqnarray*}}
\def\beq{\begin{eqnarray}}
\def\enq{\end{eqnarray}}

\usepackage[algo2e,ruled,vlined]{algorithm2e} 
 \usepackage{algorithm}
 \usepackage{algorithmicx}
\usepackage{algpseudocode}
\usepackage{ marvosym }
\usepackage{csquotes}
\usepackage{hyperref}

\def \eps{\varepsilon}

\def\argmin{\mathop{\rm argmin}}
\def\argmin_#1{\underset{#1}{\mathrm{argmin\, }}}

\def \ep{\hbox{ }\hfill$\Box$}

\newtheorem{Theorem}{Theorem}[section]

\newtheorem{Remark}[Theorem]{Remark}

\numberwithin{equation}{section} 

\begin{document}

\title{
Mean-field neural networks: learning mappings on \\ Wasserstein space
\thanks{This work is supported by  FiME, Laboratoire de Finance des March\'es de l'Energie, and the ''Finance and Sustainable Development'' EDF - CACIB Chair.}
\thanks{We thank Maximilien Germain and Mathieu Laurière for helpful discussions.}
}

\author{Huy\^en \sc{Pham}\footnote{LPSM,  Universit\'e Paris Cité,  \& FiME \sf \href{mailto:pham at lpsm.paris}{pham at lpsm.paris}} 
\and
Xavier \sc{Warin}
\footnote{EDF R\&D \& FiME \sf \href{mailto:xavier.warin at edf.fr}{xavier.warin at edf.fr}}}

\maketitle  

\begin{abstract}
We study the machine learning task  for models with  operators mapping between  the Wasserstein space of probability measures and  a space of functions,  like e.g. in  mean-field games/control problems. 
Two classes of neural networks based on  bin density and on cylindrical approximation, are proposed to learn these so-called mean-field functions,  
and are theoretically supported by universal approximation theorems.  We perform se\-veral numerical experiments for training these two mean-field neural networks, and show their accuracy and efficiency in the generalization error 
with various test distributions.  Finally, we present different  algorithms relying on mean-field neural networks for solving  time-dependent mean-field problems, and illustrate our results with numerical tests for the example of a semi-linear partial differential equation in the Wasserstein space of probability measures.  
 \end{abstract}

\section{Introduction}

Deep neural networks have been successfully used   for approximating solutions to high dimensional partial differential equations (PDEs) and control problems, and various methods either based on physics informed representation (\cite{raissi19}, \cite{sirignano2018dgm}), or probabilistic and 
backward stochastic differential equations (BSDEs)  representation (\cite{weinan2017deep}, \cite{chan2019machine}, \cite{hure2020deep}) have been recently developed in the literature, see e.g.  the  survey papers \cite{bechutjenkuc20} and \cite{GPW21}. 

\vspace{1mm}

In the last years, a novel class of control problems has emerged with the theory of mean field game/control dealing with models of large population of interacting agents. Solutions to mean-field problems are represented by functions that depend not only on the state variable of the system, but also on its probability distribution, representing the population state distribution, and can be characterized in terms of PDEs in the Wasserstein space of probability measures (called Master equation)  or BSDEs of McKean-Vlasov (MKV) type, and we refer to the two-volume monograph \cite{cardel19}, \cite{cardel2} for a comprehensive treatment of this topic. 
In such problems, the input is a probability measure on $\R^d$,  hence valued in the infinite dimensional Wasserstein space, and the output is a function defined on the support of the input probability measure.   

\vspace{1mm}

In this paper, we aim to approximate the infinite dimensional mean-field function by proposing two classes of neural network architectures. The first approach starts from the approximation of a probability measure with density by a piecewise constant density function on some given fixed 
partition of size $K$ of a truncated support of the measure, called bins, see Figure \ref{fig:binDist} in the case of a Gaussian distribution.  This allows us to approximate the infinite dimensional mapping by a function that maps an input space  of dimension $K$ corresponding to the bin density weights that can be learned by a standard deep neural network. 
We show a universal approximation theorem that justifies theoretically the use of such bin density neural network.  The second approach maps directly probability measures as input  but through a finite-dimensional neural network function in cylindrical form, for which we also state  a universal approximation theorem.  

\vspace{1mm}

Next, we show how to effectively learn mean-field function by means of these two classes of mean-field neural networks.  This is achieved by generating a data set consisting of simulated probability measures following two proposed methods, and then by training via stochastic gradient method the parameters of the mean-field neural networks. 
We perform several numerical tests for illustrating the efficiency and accuracy of these two mean-field neural networks on various examples of mean-field functions, and  we validate our results on different test distributions by computing the generalization error. 

\vspace{1mm}
    
As an application of these mean-field neural networks, we consider dynamic mean-field problems arising typically from mean-field type control, and design different algorithms of local or global type, based on regression or BSDE representation, for computing the solution. We illustrate the performance of our algorithms with 
the example of a  semi-linear PDE on the Wasserstein space.  More applications and examples from mean-field control problems and Master equations are investigated in a  companion paper \cite{phawar23b}  where we provide  a global comparison 
of the different neural network algorithms.

\vspace{1mm}

{\it Related works.} Several methods have been recently proposed for solving numerically mean field game/control problems.  We mention for instance the papers  \cite{Ruthotto9183}, \cite{linetal}  using Hamilton-Jacobi-equations and Lagrangian methods, the works by \cite{carlau19}, \cite{fouquezhang19} relying on  backward stochastic differential equations and maximum principle or the work  in \cite{gerlauphawar21a} that approximates the mean-field control problem by particle systems for reducing the problem to a finite, but possibly very  high-dimensional problem. Actually, in the latter paper, symmetry of the particle system is exploited in the numerical resolution by using a specific class of neural networks, called DeepSets \cite{deepsets}, which allows to reduce significantly the computational complexity.  However, 
in all these cited references, as the distribution probability of the state process is a deterministic function of time, the value function  and optimal control are  viewed as  functions of time and of the state, and approximated by neural networks on finite-dimensional space.  However, the solution obtained is valid for a given initial distribution of the population state, but when varying the initial distribution, the solution has to be computed again by another neural network.  In this work, we develop instead  a numerical scheme  for approximating by a suitable neural network  the solution at any initial distribution. 

\vspace{1mm}

The paper is organized as follows.  In Section \ref{sec:MF},  we formulate the learning problem,  present two network architectures: bin-density and cylindrical neural networks,  and explain the data generation and training  procedures. Numerical tests are developed in Section \ref{sec:numtest}, and applications to time dependent mean-field problems are 
given in Section \ref{sec:dynMF} with various algorithms and numerical results.  The proofs of the universal approximation theorem for mean-field neural networks are postponed in Appendix \ref{sec:appen}.

\vspace{2mm}

\noindent {\bf Notations.} Denote by $\Pc_2(\R^d)$   the Wasserstein space of square integrable probability measures equipped with the $2$-Wasserstein distance $\Wc_2$. Given $\mu$ $\in$ $\Pc_2(\R^d)$, we denote by $L^2(\mu)$ as the set of measurable functions $\phi$ on $\R^d$ s.t. 
\beqs
 |\phi|_\mu^2 &: =&  \; \int |\phi(x)|^2 \mu(\d x) \; < \; \infty. 
\enqs  
(Here $|.|$ denotes the Euclidian norm). 
Given some  $\mu$ $\in$ $\Pc_2(\R^d)$, and $\phi$ a measurable function on $\R^d$ with quadratic growth condition, hence in $L^2(\mu)$, we set: $\E_{X\sim\mu}[\phi(X)]$ $:=$ $\int \phi(x) \mu(\d x)$.  We also denote by $\bar\mu$ $:=$ $\E_{X\sim\mu}[X]$.

\section{Learning mean-field functions}  \label{sec:MF}

Given a function $V$ on $\R^d\times\Pc_2(\R^d)$, valued on $\R^p$, with quadratic growth condition w.r.t. the first argument in $\R^d$, we aim to approximate the infinite-dimensional mapping
\begin{align} \label{defVc} 
\Vc : \mu \in \Pc_2(\R^d) & \longmapsto \;  V(\cdot,\mu) \in L^2(\mu), 
\end{align} 
called mean-field function, by a map  $\Nc$ constructed from suitable classes of neural networks. The mean-field network $\Nc$ takes inputs composed of two parts: $\mu$ a probability measure and $x$ in the support of $\mu$, and outputs $\Nc(\mu)(x)$. 
The quality of this approximation is measured by the error:
\beqs
L(\Nc) & := & \int_{\Pc_2(\R^d)}  \Ec_{\Nc}(\mu)   \nu(\d \mu),  \\
\;  \mbox{ with } \quad   \Ec_{\Nc}(\mu) & := &  \big| \Vc(\mu) - \Nc(\mu) \big|_\mu^2 \;  = \; \E_{X\sim\mu} \big| V(X,\mu) - \Nc(\mu)(X) \big|^2, 
\enqs
where $\nu$ is a probability measure on $\Pc_2(\R^d)$, called training measure. The learning of  the mean-field functional $\Vc$ will be then performed by minimizing over the parameters of the neural network operator $\Nc$ the loss function
\begin{align} \label{lossM}
L_M(\Nc) \; := \; \frac{1}{M} \sum_{m=1}^M   \Ec_{\Nc}(\mu^{(m)}), 
\end{align}
where $\mu^{(m)}$, $m$ $=$ $1,\ldots,M$ are training samples of $\nu$.  We denote by $\widehat\Nc^M$ the learned functional from this minimization problem,  and for test data $\mu^{test}$ (different from the training data set $(\mu^{(m)})_m$), we shall compute the test (generalization) error $\Ec_{\widehat\Nc^M}(\mu^{test})$.

 \subsection{Neural networks approximations} \label{secNN}

\begin{figure}[H]
    \centering
    \includegraphics[width=0.4\textwidth]{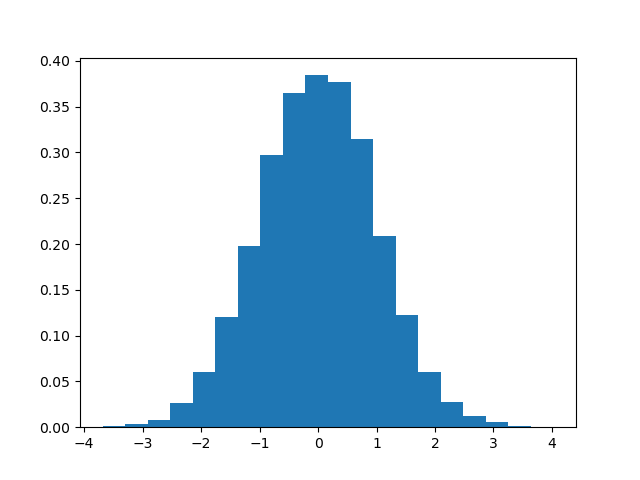}
    \caption{Bin approximation of a Gaussian distribution.}
    \label{fig:binDist}
\end{figure}

\paragraph{{\bf Bin density-based approximation}}  Let us denote by $\Dc_2(\R^d)$ the subset of probability measures $\mu$ in $\Pc_2(\R^d)$ which admit density functions $\mrp^\mu$  with respect to the Lebesgue measure $\lambda_d$ on $\R^d$. 
Fix  $\Kc$ as  a bounded rectangular domain in $\R^d$, and divide $\Kc$ into a number $K$ of bins, Bin$(k)$, $k$ $=$ $1,\ldots,K$: $\cup_{k=1}^K {\rm Bin}(k)$ $=$ $\Kc$, of center $x_k$, and  with same area size $h$ $=$ $\lambda_d(\Kc)/K$.   
Given $\mu$ $\in$ $\Dc_2(\R^d)$, we consider the bin approximation of its density function (see figure \ref{fig:binDist}), that is the truncated piecewise-constant density function defined on $\Kc$ by  
\beqs
\hat\mrp^\mu_{\Kc}(x) &=& p^\mu_k :=  \frac{\mrp^\mu(x_k)}{\sum_{k=1}^K \mrp^\mu(x_k) h}, \;  \mbox{ if } x \in {\rm Bin}(k), \; k = 1,\ldots,K, 
\enqs
and $\hat\mrp^\mu_{\Kc}(x)$ $=$ $0$ for $x \in \R^d \setminus\Kc$,  set  $\bop^\mu$ $:=$ $(p_k^\mu)_{k\in\llbracket 1,K\rrbracket}$, which lies in $\Dc_K$ $:=$ $\{ \bop = (p_k)_{k\in\llbracket 1,K\rrbracket} \in \R_+^K: \sum_{k=1}^K p_k h =1\}$, and 
called  density bins of  the probability measure in $\Dc_2(\R^d)$  of density function $\hat\mrp_\Kc^\mu$, denoted by $\hat\mu^K$ with support on $\Kc$ 

Conversely, given $\bop$ $=$ $(p_k)_{k\in\llbracket 1,K\rrbracket}$ $\in$ $\Dc_K$,  we can associate the piecewise-constant density function defined on $\R^d$ by 
\begin{align} \label{defdis} 
\mrp(x) &=\;  p_k, \;  \mbox{ if } x \in {\rm Bin}(k), \; k = 1,\ldots,K, \quad \mrp(x) = 0, \quad x \in \R^d \setminus\Kc. 
\end{align}
We then denote by $\mu$ $=$ $\Lc_D\big(\bop  \big)$  the bin density  probability measure on $\Pc_2(\R^d)$ with piecewise-constant density function $\mrp$ as in \eqref{defdis}, hence with support on $\Kc$, and we note that $\hat\mu^K$ $=$  $\Lc_D(\bop^\mu)$.


 
 A mean-field density-based network is an operator on $\Dc_2(\R^d)$ in the form
 \beqs
 \Nc_D(\mu) &=&  \Phi(\cdot,\bop^\mu),
 \enqs
 where  $\Phi$ $=$ $\Phi_\theta$ is a neural network function from $\R^d\times\Dc_K$ into $\R^p$, whose architecture can be constructed as follows: 
\begin{enumerate}
\item[(i)]  Classical feedforward neural network, i.e. in the form:
\beqs
(x,\bop) \in \R^d\times\Dc_K & \mapsto & \Phi_\theta(x,\bop) \; = \;  \Ac_{L+1} \circ  \underbrace{ \sigma \circ  \Ac_{L} }_{L-\mbox{layer}}  \circ  \ldots  \circ 
\underbrace{\sigma \circ  \Ac_1(x,\bop)}_{1-layer} \; \in \; \R^{p},  \\
\Ac_\ell (x,\bop)  &=&  w_\ell  \left( \begin{array}{c} 
                                                     x \\
                                                     \bop 
                                                     \end{array} \right)  + b_\ell \; \in \; \R^{d_\ell}, \quad d_{L+1} = p, 
\enqs  
with $L$ hidden layers (layer $\ell$ with $d_\ell$ neurons), parameters $\theta$ $=$ $(w_\ell,b_\ell)_{\ell}$: $w_\ell$ weight, $b_\ell$  bias, 
an activation function $\sigma$ from $\R$ into $\R$ (composition is componentwise), like e.g.  $\tanh$, sigmoid, or Relu. 
\item[(ii)] DeepOnet structure (see \cite{lu2019deeponet}):  $\Phi_\theta(x,\bop)$ $=$ $\sum_{\ell=1}^L \mrb_\ell \mrt_\ell$, where $(\mrb_\ell)_\ell$ is the ouput of a branch net with input $\bop$ $=$ $(p_k)_k$ representing the sensors, and $(\mrt_\ell)_\ell$ is  the output of a trunk net with input $x$.  
\item[(iii)] Other structures like the networks developed in  \cite{warin2021reservoir} for a differentiated treatment of uncertainties and storage level in reservoir optimization.
\end{enumerate}

 Let us denote by $\Dc_{c}(\Kc)$ the subset of elements $\mu$ in $\Dc_2(\R^d)$ with support in $\Kc$, with continuous density functions $\mrp^\mu$, and for $\mu$ $\in$ $\Dc_{c}(\Kc)$, we set $\omega_\Kc^\mu$ as  its modulus of uniform continuity on $\Kc$. 
 Given a modulus of continuity $\bar\omega$, i.e. a nondecreasing function on $\R_+$  s.t.  $\lim_{t\downarrow 0} \bar\omega(t)$ $=$ $\bar\omega(0)$ $=$ $0$, we denote by $\Dc_c^{\bar\omega}(\Kc)$ the subset of elements $\mu$ $\in$ $\Dc_c(\Kc)$ such that 
 $\omega_\Kc^\mu$ $\leq$ $\bar\omega$ on a neighborhood of $t$ $=$ $0$.

The justification for the use of  bin-density neural networks is due to the following universal approximation theorem. 

\begin{Theorem} \label{theounivbin} 
Let $\Kc$ be bounded rectangular domain in  $\R^d$, $\bar\omega$ a modulus of continuity, and 
$V$  a continuous function on $\R^d\times\Pc_2(\R^d)$. Then, for all 
$\varepsilon$ $>$ $0$, there exists $K$ $\in$ $\N^*$, and 
$\Phi$ a neural network on $\R^d\times\R^K$ such that 
\begin{align}
\big| V(x,\mu) - \Phi(x,\bop^\mu)\big| 
& \leq  \; \varepsilon, \quad \forall x \in \Kc, \; \mu \in . 
\end{align}
\end{Theorem}
\noindent  {\bf Proof.}
 See Appendix \ref{sec:appen}. 
 \ep

 \paragraph{Cylindrical  approximation.}  
 A mean-field cylindrical neural network is an operator on $\Pc_2(\R^d)$ in the form
 \begin{align}
\Nc_C(\mu) &= \;  \Psi_\theta(\cdot, < \varphi_\theta,\mu>),
 \end{align}
 where $\Psi_\theta$ is a feedforward network function from $\R^d\times\R^k$ into $\R^p$, and $\varphi_\theta$ is another feedforward network function from $\R^d$ into $\R^k$ (called latent space).  Here we denote $<\varphi_\theta,\mu>$ $:=$ $\int \varphi_\theta(x) \mu(\d x)$ $=$ $\E_{X\sim\mu}[\varphi_\theta(X)]$.  
 By misuse of language, we call $(\Psi_\theta,\varphi_\theta)$ such cylindrical neural network with $\varphi_\theta$ the inner network, and $\Psi_\theta$ the outer network.

\vspace{1mm}

We state  a universal approximation theorem  that justifies the use of cylindrical mean-field neural networks.  It is stated with an $L^2$-distance, which is, in practice, the distance that is minimized during the training process. 

\begin{Theorem} \label{theounivcyl2} 
Let $\nu$ be a probability measure on $\Pc_2(\R^d)$, and $V$ be a continuous function from $\R^d\times\Pc_2(\R^d)$ into $\R^p$ s.t. $\|V\|^2_{L^2(\nu)}$ $:=$ 
$\int_{\Pc_2(\R^d)} |V(.,\mu)|_\mu^2 \nu(\d \mu)$ $<$ $\infty$.  Then, for all $\varepsilon$ $>$ $0$, there exists 
$k$ $\in$ $\N$, $\Psi$ a  neural network from $\R^d\times\R^k$ into $\R^p$, $\varphi$ a neural network from $\R^d$ into $\R^k$ such that
\begin{align}
\int_{\Pc_2(\R^d)}  
 \big| V(.,\mu) - \Psi(.,<\varphi,\mu>) \big|_\mu^2 \nu(\d \mu) 
& \leq \; \varepsilon. 
\end{align}
\end{Theorem}
\noindent  {\bf Proof.}
See Appendix \ref{sec:appen}. 
\ep

 \subsection{Data generation}
 \label{subsec:Data}
 
 The training of neural networks for approximating mean-field function  relies on samples of probability measures $\mu$ and of random state value $X$ distributed according to $\mu$. We propose  two  methods.  
\begin{enumerate}
\item We draw a random grid $\bolx$ $=$ $(x_k)_{k\in\llbracket 1,K\rrbracket}$  of $K$ points in $\R^d$, according to some probability measure on $\R^d$, e.g. uniform distribution on $\Kc$, as well as a random element $\bopi$ $=$ $(\pi_k)_{k\in\llbracket 1,K\rrbracket}$ in  the simplex $\Sc_K$ $=$ 
$\{ \bopi = (\pi_k)_{k\in\llbracket 1,K\rrbracket} \in \R_+^K: \sum_{k=1}^K \pi_k  =1\}$.  This can be 
done for example from a sample $e_1,\ldots,e_K$ of positive random variables  according to an exponential law, and by setting $ \pi_k$ $=$ $\frac{e_k}{\sum_{k=1}^K e_k}$, $k$ $=$  $1,\ldots,K$. This generates a (random) quantized  probability measure:
\beqs
\Lc_Q^{\bolx}(\bopi) & := & \sum_{k=1}^K \pi_k \delta_{x_k}, 
\enqs
that is the discrete probability measure with support on the grid $\bolx$ and with probability weights $\bopi$. 
\item We draw random vector $\bop$ $=$ $(p_k)_{k\in\llbracket 1,K\rrbracket}$ in $\Dc_K$.   This can be  done for example from a sample $e_1,\ldots,e_K$ of positive random variables  according to an exponential law on $\R_+$, e.g. with parameter 1, and by setting 
$p_k$ $=$ $\frac{e_k}{\sum_{k=1}^K e_k h}$, $k$ $=$  $1,\ldots,K$.  This generates (random) bin density probability measure $\mu$ $=$ $\Lc_D(\bop)$,  whose cumulative distribution function is given in the one-dimensional case ($d$ $=$ $1$, ${\rm Bin}(k)$ $=$ $[x_{k-1},x_k)$, $k$ $=$ $1,\ldots,K$) by
\begin{equation}
F_{\bop}(x) = \begin{cases}
0,  &\quad x < x_0 \\
\Sum_{j=1}^{k-1} p_j h_j \; + \; p_k(x-x_{k-1}), & \quad x \in {\rm Bin}(k), \; k=1,\ldots,K, \\
1,&  \quad x \geq  x_K, 
\end{cases}
\end{equation}
 with the convention that $\sum_{j=1}^{k-1}$ $=$ $0$ for $k$ $=$ $1$.  Its inverse function is then explicitly given by 
 \beqs
 F_{\bop}^{-1}(u) &=& x_{k_u-1} + \frac{ u - \Sum_{j=1}^{k_u-1} p_j h_j}{p_{k_u}},  \quad u \in [0,1],  \\
 \mbox{ with }\quad  \quad  k_u  &:=&  \inf\{ k \in \llbracket 1,K \rrbracket:  \sum_{j=1}^k p_j h_j \geq u \}. 
 \enqs
 We then draw an  uniform random variable $U$ on $[0,1]$, which generates a random variable $X$ $=$ $ F_{\bop}^{-1}(U)$  distributed  according to  $\Lc_D(\bop)$. 
\end{enumerate}
 There are other methods for generating random probability measures that  can be found in the literature on probability theory, and we refer e.g. to an overview in \cite{mon01}, but it is not clear how one can simulate easily random variables distributed according to such random probability measures.

\subsection{Training mean-field neural networks}  \label{sec:trainNN}


According to the choice of the mean-field neural network,   the training for learning the mean-field operator $\Vc$ in \eqref{defVc} is performed as follows:
\begin{enumerate}

\item {\it Bin density-based neural network}. We draw a sample $\bop^{(m)}$  of vectors in $\Dc_K$, which gene\-rates a sample of bin density probability measures $\mu^{(m)}$ $=$ $\Lc_D(\bop^{(m)})$, $m=1,\ldots,M$. 
By noting that the  density bins of the density probability measure $\mu^{(m)}$ is $\bop^{(m)}$, the training of the bin density-based  neural network $\Nc_D$ from the minimization of the loss function in \eqref{lossM} consists in 
minimizing over the parameters $\theta$ of a feedforward neural network $\Phi_\theta$ on $\R^d\times\Dc_K$ the loss function  
\beqs
L_D(\theta) \; = \; \frac{1}{M} \sum_{m=1}^M  \E_{X\sim\mu^{(m)}} \big| V(X,\mu^{(m)}) -  \Phi_\theta(X,\bop^{(m)}) \big|^2. 
\enqs
For the approximation of this expectation $\E_{X\sim\mu^{(m)}}[.]$  when applying  SGD we shall use for each $m$, a batch $X^{(n)}$, $n$ $=$ $1,\ldots,N$, of samples of $X$ $\sim$ $\mu^{(m)}$. 
Notice that this method works effectively in dimension $d$ $=$ $1$  in order to be able to simulate $X$. 
\item {\it Cylindrical neural network}. We draw a sample $\mu^{(m)}$, $m=1,\ldots,M$, of proba\-bility measures in $\Pc_2(\R^d)$, either according to discrete probability measures $\mu^{(m)}$ $=$ $\Lc_Q^{\bolx^{(m)}}(\bopi^{(m)})$,
 or to bin density probability measures  $\mu^{(m)}$ $=$ $\Lc_D(\bop^{(m)})$, and minimize over the parameters $\theta$ of a cylindrical neural network $(\Psi_\theta,\varphi_\theta)$ the loss function: 
\beqs
L_C(\theta) \; = \; \frac{1}{M} \sum_{m=1}^M \E_{X\sim\mu^{(m)}}  \Big| V(X,\mu^{(m)}) - \Psi_\theta\big(X, \E_{X\sim\mu^{(m)}}[ \varphi_\theta(X)] \big) \Big|^2.
\enqs
Again, when applying  SGD for the approximation of this expectation $\E_{X\sim\mu^{(m)}}[.]$, we shall use for each $m$, a batch $X^{(n)}$, $n$ $=$ $1,\ldots,N$, of samples of $X$ $\sim$ $\mu^{(m)}$. 
\end{enumerate}

\section{Numerical tests} \label{sec:numtest}

We test our two choices of mean-field neural networks by computing the corresponding training mean square error (MSE) and test (generalization) error for different cases of mean-field functions $V$ on $\R\times\Pc_2(\R)$.

We first consider the two  following cases of mean-field functions: 
\begin{itemize}
\item[A.] {\it Case A: a quadratic function of the measure}  
    \begin{align}
    V(x, \mu) & = \;  x + \bar \mu + 2 {\rm Var}(\mu),  
    \end{align}
where $\bar\mu$ $:=$ $\E_{X\sim\mu}[X]$, ${\rm Var}(\mu)$ $:=$ $\E_{X\sim\mu}[X^2]- |\bar\mu|^2$. 
    \item[B.] {\it Case B : a superquantile  function}
    \begin{align}
        V(x, \mu)  &= \;  x^2+ 2 x \E_{X \sim \mu}[ X | X \ge Q_{\mu}(0.5)] + \E_{X \sim \mu}[ X | X \ge Q_{\mu}(0.5)]^2, 
    \end{align}
    where $Q_\mu(p) = \inf\{x \in \R : p \le F_\mu(x) \}$ and $F_{\mu}$ is the cumulative distribution function  of a random variable with law  $\mu$.
    \end{itemize}

We first  illustrate  the convergence of the bin method with classical feedforward neural network for different  hyperparameters. The data are generated in all cases with the second method described in section \ref{subsec:Data} for $K$ $=$ $100$ bins.  We take $\hat M=20$ distributions as the batch size during training with $N$ $=$ $50000$, and 
give  the convergence of the ADAM methods \cite{kingma2014adam} by plotting the moving average on 10 values   of the  MSE calculated  every $100$ iterations of SGD  with $M=1000$ testing distributions. The initial learning rate is $10^{-3}$ and we use tensorflow \cite{2015tensorflow}.
 \begin{figure}[H]
     \centering
     \begin{minipage}[t]{0.40\linewidth}
  \centering
 \includegraphics[width=\textwidth]{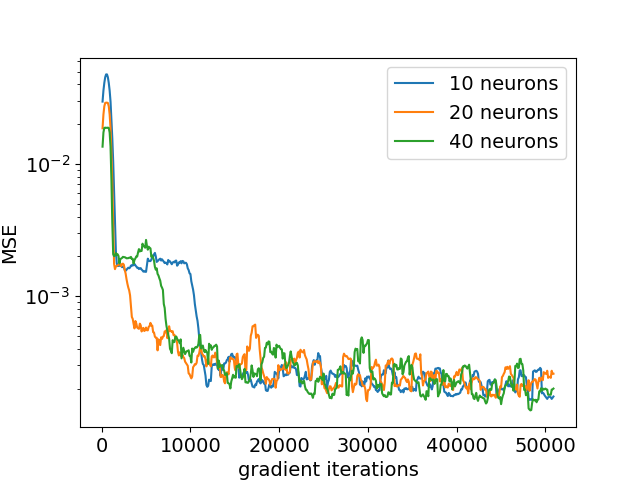}
\caption*{ReLu, 2 hidden layers}
\end{minipage}
 \begin{minipage}[t]{0.40\linewidth}
  \centering
 \includegraphics[width=\textwidth]{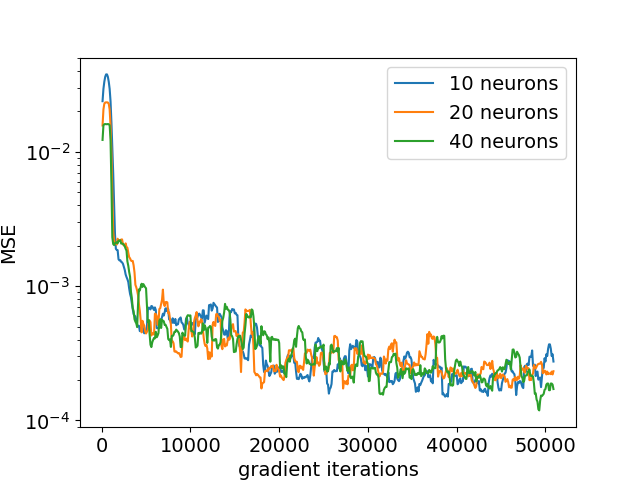}
\caption*{ReLu, 3 hidden layers}
\end{minipage}
   \centering
     \begin{minipage}[t]{0.40\linewidth}
  \centering
 \includegraphics[width=\textwidth]{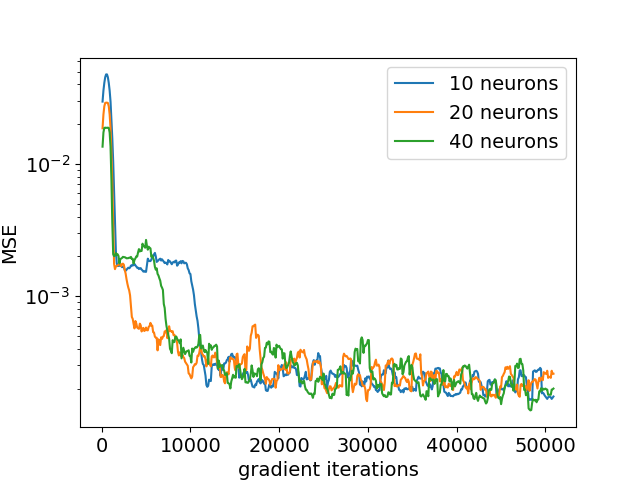}
\caption*{Tanh, 2 hidden layers}
\end{minipage}
 \begin{minipage}[t]{0.40\linewidth}
  \centering
 \includegraphics[width=\textwidth]{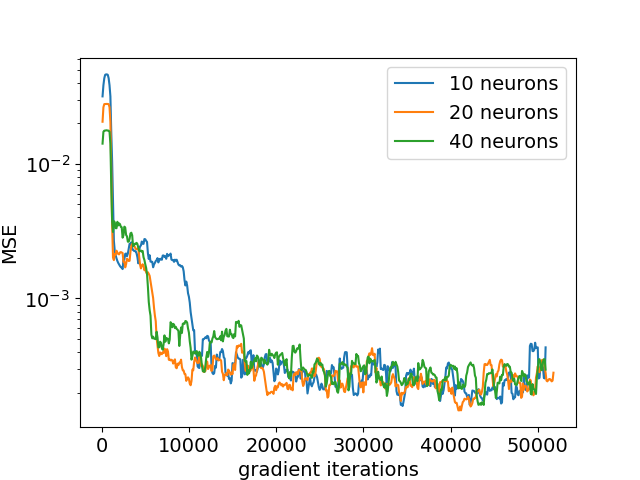}
\caption*{Tanh, 3 hidden layers}
\end{minipage}
     \caption{\footnotesize{Bin approximation: Training error  for case A depending on the number of neurons.}}
     \label{fig:binConverA}
 \end{figure}

 \begin{figure}[H]
     \centering
     \begin{minipage}[t]{0.40\linewidth}
  \centering
 \includegraphics[width=\textwidth]{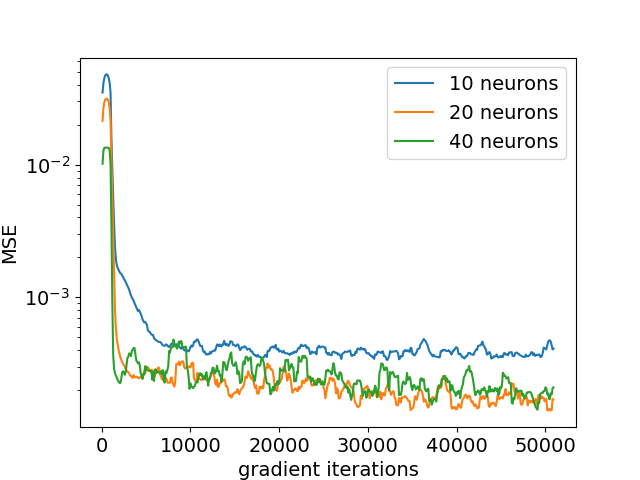}
\caption*{ReLu, 2 hidden layers}
\end{minipage}
 \begin{minipage}[t]{0.4\linewidth}
  \centering
 \includegraphics[width=\textwidth]{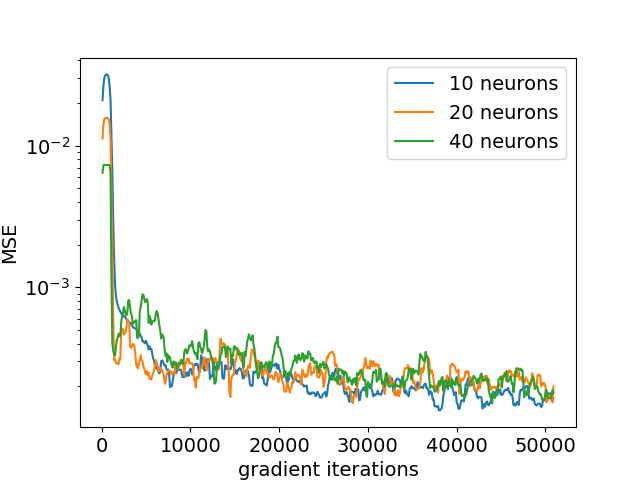}
\caption*{ReLu, 3 hidden layers}
\end{minipage}
   \centering
     \begin{minipage}[t]{0.4\linewidth}
  \centering
 \includegraphics[width=\textwidth]{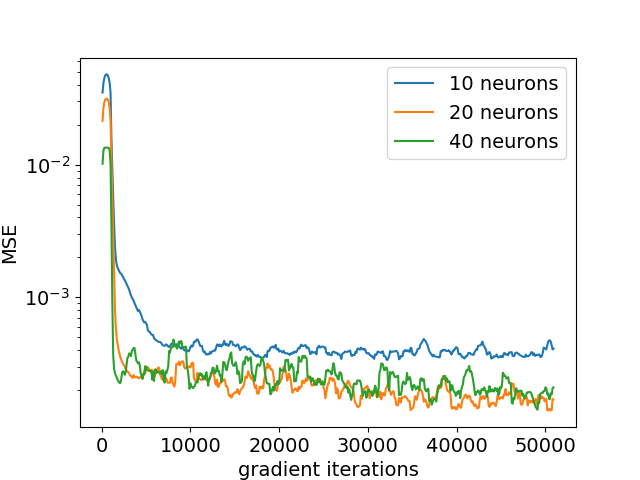}
\caption*{Tanh, 2 hidden layers}
\end{minipage}
 \begin{minipage}[t]{0.4\linewidth}
  \centering
 \includegraphics[width=\textwidth]{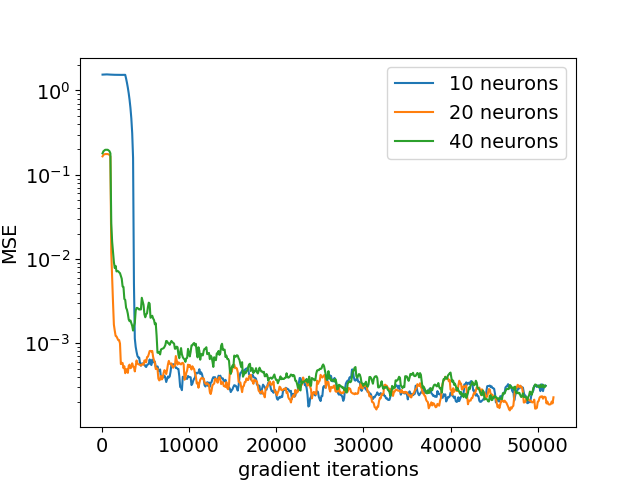}
\caption*{Tanh, 3 hidden layers}
\end{minipage}
     \caption{\footnotesize{Bin approximation: Training error  for case B depending on the number of neurons.}}
     \label{fig:binConverB}
 \end{figure}

\begin{figure}[H]
     \centering
     \begin{minipage}[t]{0.4\linewidth}
  \centering
 \includegraphics[width=\textwidth]{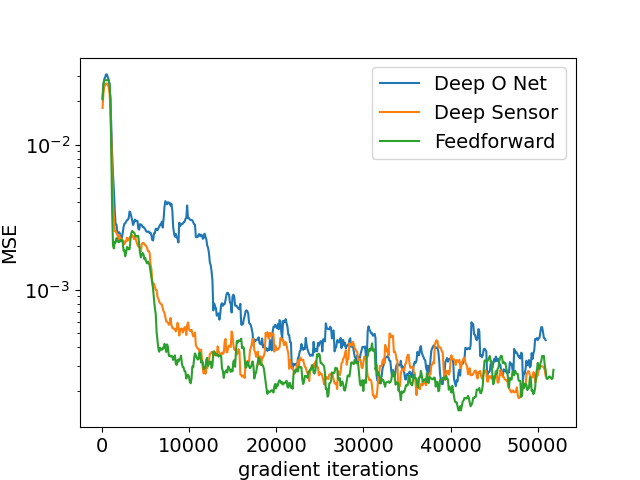}
\caption*{Case A}
\end{minipage}
 \centering
     \begin{minipage}[t]{0.4\linewidth}
  \centering
 \includegraphics[width=\textwidth]{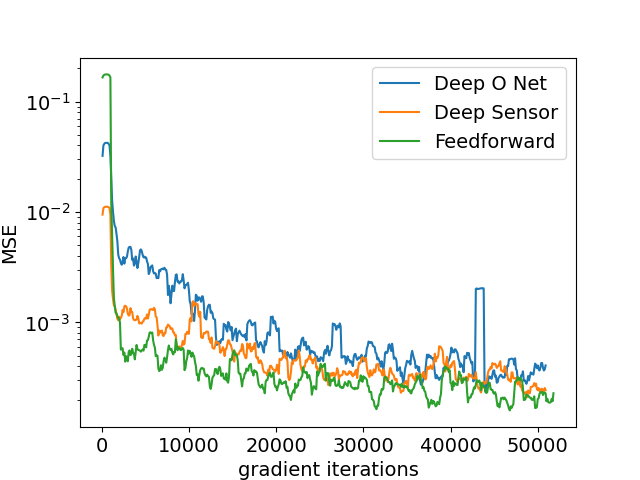}
\caption*{Case B}
\end{minipage}
\caption{\label{fig:testNetForBin} Comparison of different networks for bin approximation.}
\end{figure}

 
 Results for the bin approximation on Figures \ref{fig:binConverA}, \ref{fig:binConverB} 
 indicate that three layers with 20 neurons and a tanh activation function is a good choice for these regular functions.
We may also wonder if another architecture for neural network  may improve the results: we test the DeepONet network and the network developed in paragraph 3.2 in
 \cite{warin2021reservoir} (Deep Sensor in the graphs).  Results with three hidden layers with 20 neurons for each network, the tanh activation function,  are given on Figure \ref{fig:testNetForBin}, and show that other networks do not  seem to improve the feedforward results.  In the sequel we only use the classical feedforward network.


\begin{figure}[H]
     \centering
     \begin{minipage}[t]{0.35\linewidth}
  \centering
 \includegraphics[width=\textwidth]{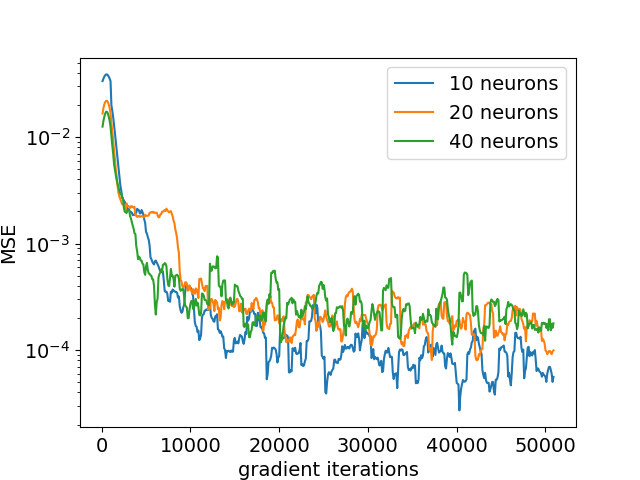}
\caption*{ReLu, $L=2$ , $Q=10$}
\end{minipage}
 \begin{minipage}[t]{0.35\linewidth}
  \centering
 \includegraphics[width=\textwidth]{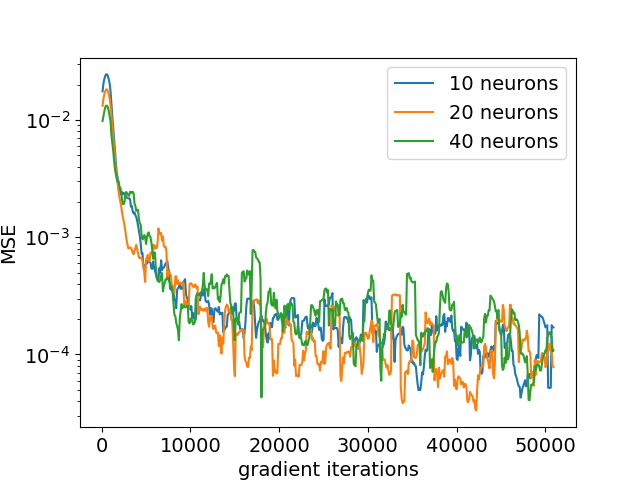}
\caption*{ReLu, $L=2$ , $Q=20$}
\end{minipage}
     \begin{minipage}[t]{0.35\linewidth}
  \centering
 \includegraphics[width=\textwidth]{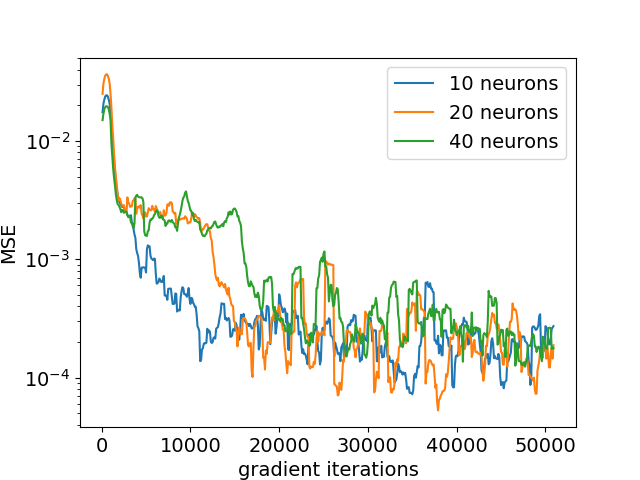}
\caption*{ReLu, $L=3$ , $Q=10$}
\end{minipage}
 \begin{minipage}[t]{0.35\linewidth}
  \centering
 \includegraphics[width=\textwidth]{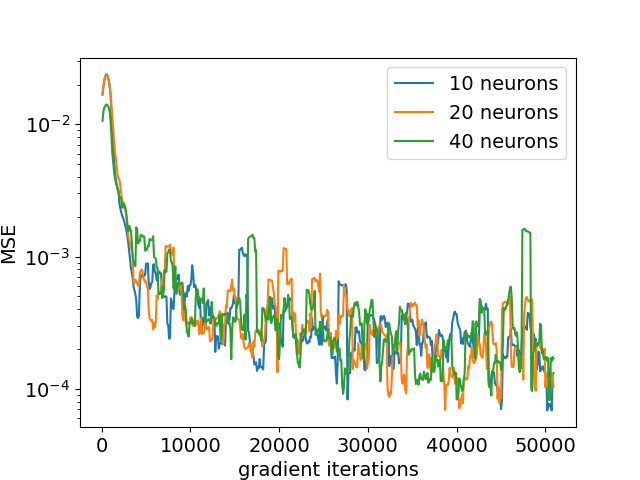}
\caption*{ReLu, $L=3$ , $Q=20$}
\end{minipage}
 \begin{minipage}[t]{0.35\linewidth}
  \centering
 \includegraphics[width=\textwidth]{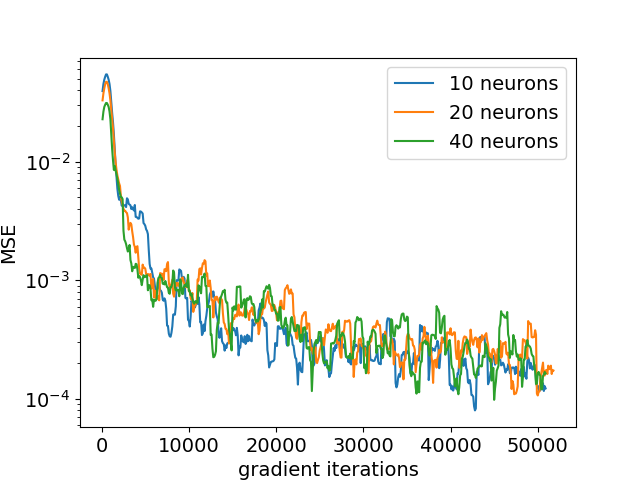}
\caption*{Tanh, $L=2$ , $Q=10$}
\end{minipage}
 \begin{minipage}[t]{0.35\linewidth}
  \centering
 \includegraphics[width=\textwidth]{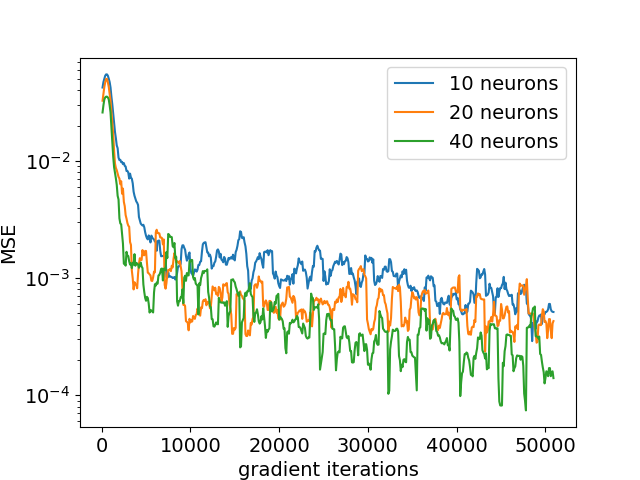}
\caption*{Tanh, $L=2$ , $Q=20$}
\end{minipage}
     \begin{minipage}[t]{0.35\linewidth}
  \centering
 \includegraphics[width=\textwidth]{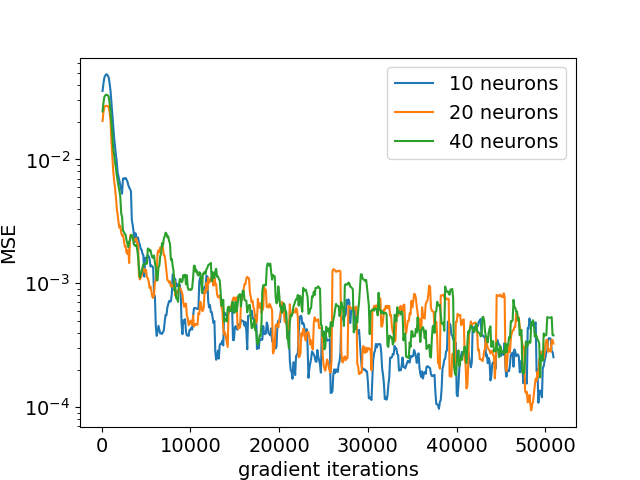}
\caption*{Tanh, $L=3$ , $Q=10$}
\end{minipage}
 \begin{minipage}[t]{0.35\linewidth}
  \centering
 \includegraphics[width=\textwidth]{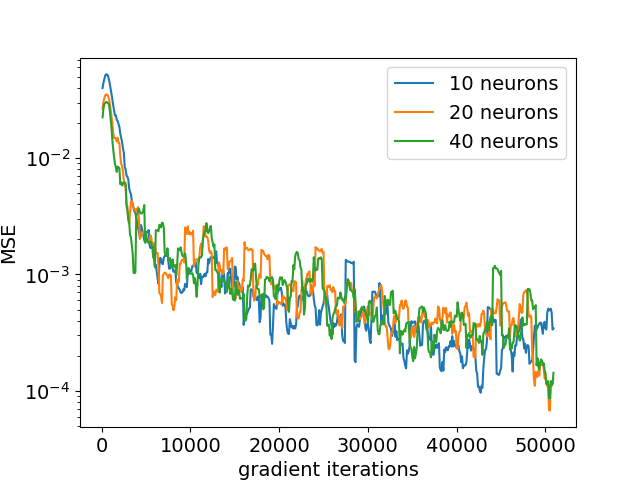}
\caption*{Tanh, $L=3$ , $Q=20$}
\end{minipage}
     \caption{Cylinder network: Training error  for case A depending on the number of neurons $k$ $=$ $10$, $20$ or $40$ of the inner network.}
     \label{fig:cylinderConverB}
 \end{figure}

Next, we consider the cylindrical approximation, which uses two networks:  
One inner network $\varphi$ having $L$ layers of $k$ neurons (so with output in dimension $k$),  and one outer network $\Psi$  having $L$  hidden  layers with $Q$ neurons. 
 The convergence of the training error is illustrated in Figure  \ref{fig:cylinderConverB} for the case A: it  indicates that a tanh activation function using 2 layers, $Q=10$, $k=20$ is a good choice. This result is confirmed  on test case B but not reported here. 
 

 In the sequel, all results for case A and B are obtained using the previously fitted networks.  
 We have shown that the global training  errors decrease correctly, and we shall now compute  the error for various given test distributions, and for the two choices of mean-field neural networks: 
 \begin{enumerate}
 \item {\it Bin density-based neural network}. We estimate  the density bins $\bop^{test}$ of $\mu^{test}$ from samples  $X^{(n)}$, $n$ $=$ $1,\ldots,N$, of $\mu^{test}$ as
 \beqs
 p_k^{test} &=&  \frac{ \#  \{ n \in \llbracket 1,N\rrbracket:  {\rm Proj}_\Kc(X_{}^{(n)}) \in \mbox{Bin}(k) \} }{Nh}, \quad k=1,\ldots,K,
 \enqs
 where ${\rm Proj}_{\Kc}(.)$ is the projection on $\Kc$.  We then compute the error test
 \beqs
 \Ec_{\widehat{\Nc}^M}(\mu^{test}) &=& \E_{X \sim \mu^{test}}  \big|  V(X,\mu^{test}) - \Phi_{\hat\theta_M}(X,\bop^{test}) \big|^2  \nonumber \\
 & \simeq & \frac{1}{N} \sum_{n=1}^N  \big|  V(X^{(n)},\mu^{test}) - \Phi_{\hat\theta_M}(X^{(n)},\bop^{test}) \big|^2. 
 \label{eq:errBins}
 \enqs
\item {\it Cylindrical neural network}.  From  samples  $X^{(n)}$, $n$ $=$ $1,\ldots,N$, of $\mu^{test}$, we shall next compute the error test as
\begin{footnotesize}
\beqs
 \Ec_{\widehat{\Nc}^M}(\mu^{test}) &=&  \E_{X \sim \mu^{test}}  \big|  V(X,\mu^{test}) -   \Psi_{\hat\theta_M}\big(X^{(n)}, \E_{X\sim\mu^{test}}[\varphi_{\hat\theta_M}(X)]\big)   \big|^2  \nonumber\\
 & \simeq &  \frac{1}{N} \sum_{n=1}^N  \big| V(X^{(n)},\mu^{test}) -    \Psi_{\hat\theta_M}\big(X^{(n)}, \frac{1}{N} \sum_{i=1}^N \varphi_{\hat\theta_M}(X^{(i)}) \big) \big|^2. 
\enqs  
\end{footnotesize}
 \end{enumerate}
 
\begin{figure}[H]
\begin{minipage}[t]{0.32\linewidth}
  \centering
 \includegraphics[width=\textwidth]{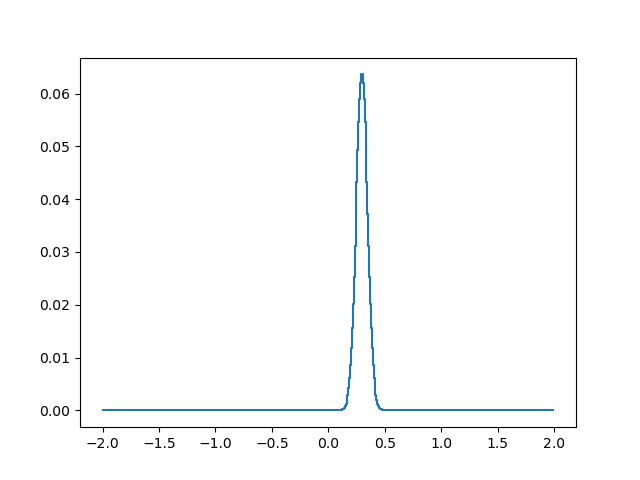}
\caption*{Test 1}
\end{minipage}
\begin{minipage}[t]{0.32\linewidth}
  \centering
 \includegraphics[width=\textwidth]{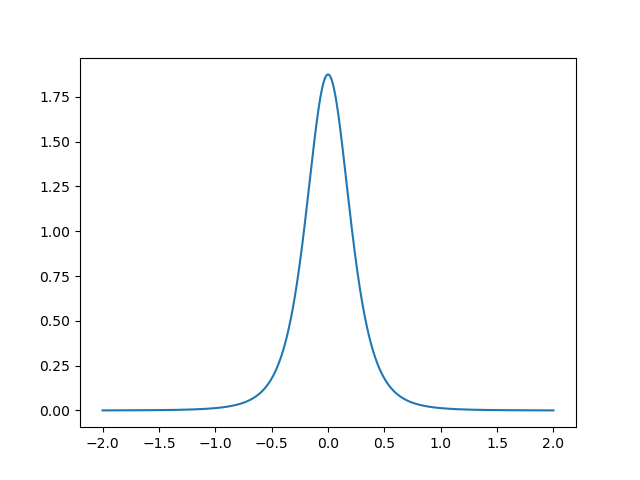}
\caption*{ Test 2}
\end{minipage}
\begin{minipage}[t]{0.32\linewidth}
  \centering
 \includegraphics[width=\textwidth]{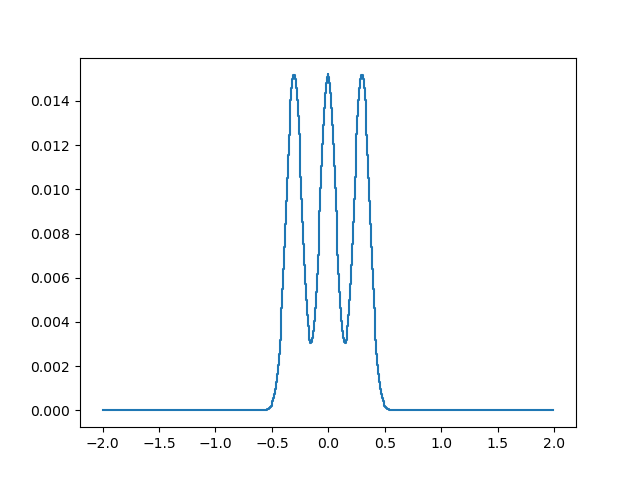}
\caption*{ Test 3}
\end{minipage}
\caption{\label{fig:distTest} Distributions used to test approximation algorithms.}
\end{figure}

We test three distributions of $X^{test}$ $\sim$  $\mu^{test}$ plotted in Figure \ref{fig:distTest} and given by: 
\begin{itemize}
    \item[(i)]  Test 1 : Gaussian with $ \bar \mu^{test}  =0.3$, $std(\mu^{test})=0.05$. 
    \item[(ii)]  Test 2 : A student distribution  with parameter $\nu=4$ and  localization factor $0$, scale factor $\sigma= 0.2$, therefore with a variance $\frac{\nu}{\nu-2} \sigma^2$.
    \item[(iii)]  Test  3 : Mixture of three  gaussians: $X_0 =  a[ - 1_{\lfloor 3U \rfloor= 0}   +  1_{\lfloor 3U \rfloor= 1 }  ] +  b Y $ with $U \sim {\cal U}(0,1)$, $a=0.3$, $b = 0.07$, $Y \sim  \mathcal{N}(0,1)$ independent of $U$.  
\end{itemize}

\begin{figure}[H]
     \centering
     \begin{minipage}[t]{0.325\linewidth}
  \centering
 \includegraphics[width=\textwidth]{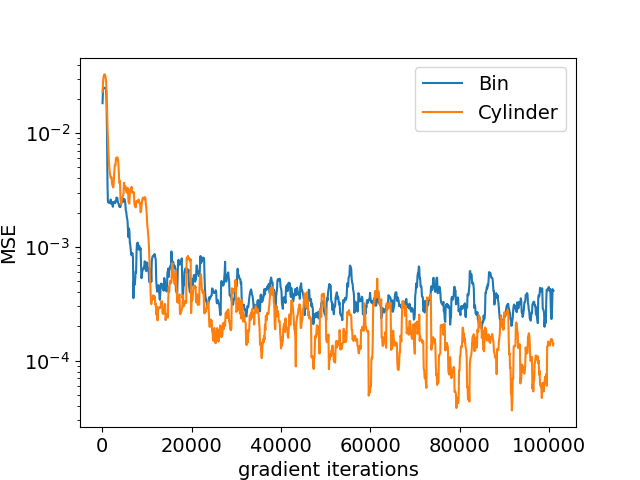}
\caption*{\tiny Training MSE : N=10000}
\end{minipage}
\begin{minipage}[t]{0.325\linewidth}
  \centering
 \includegraphics[width=\textwidth]{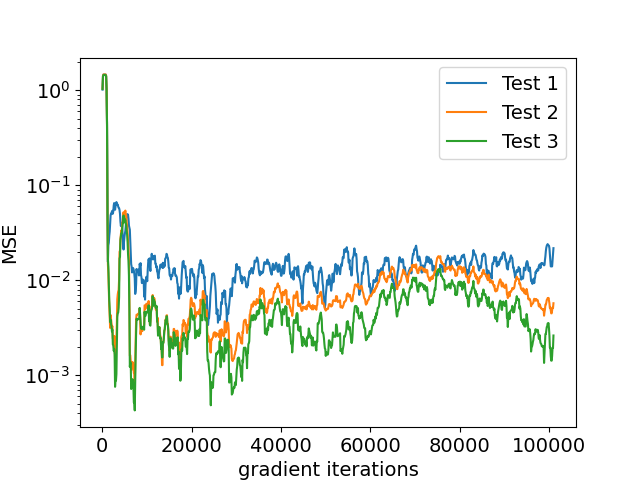}
\caption*{\tiny Test error with   bins : N=10000}
\end{minipage}
  \begin{minipage}[t]{0.325\linewidth}
  \centering
 \includegraphics[width=\textwidth]{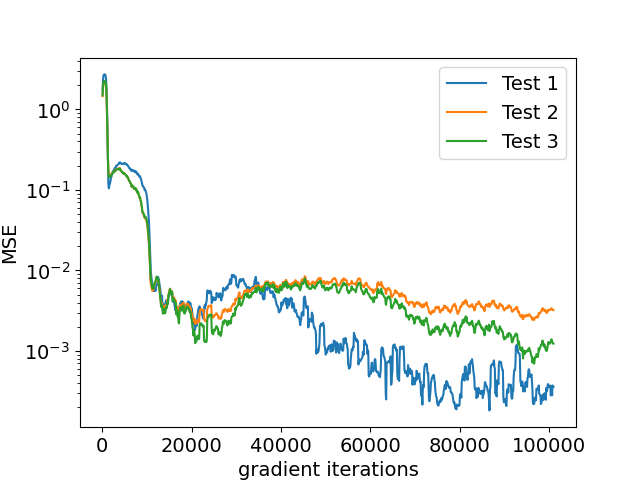}
\caption*{\tiny Test error with cylinder: N=10000}
\end{minipage}
     \begin{minipage}[t]{0.325\linewidth}
  \centering
 \includegraphics[width=\textwidth]{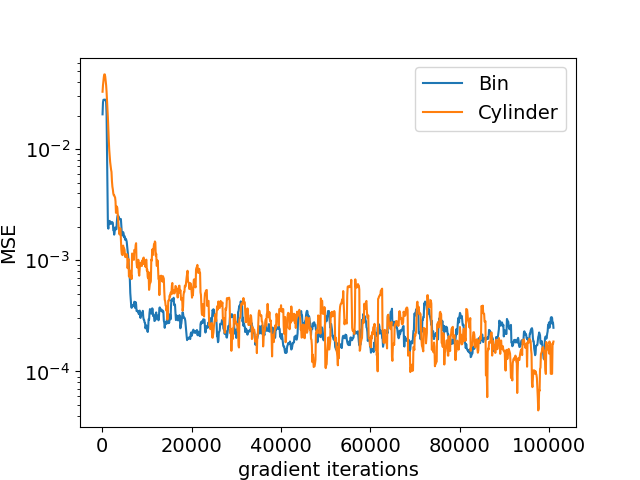}
\caption*{\tiny Training MSE : N=50000}
\end{minipage}
\begin{minipage}[t]{0.325\linewidth}
  \centering
 \includegraphics[width=\textwidth]{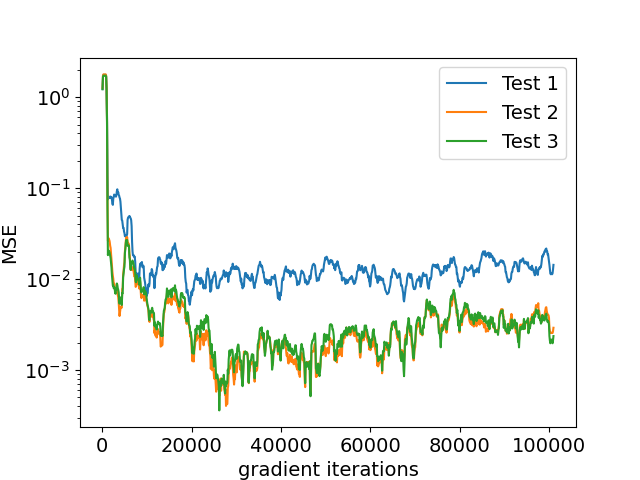}
\caption*{ \tiny Test error with   bins : N=50000}
\end{minipage}
  \begin{minipage}[t]{0.325\linewidth}
  \centering
 \includegraphics[width=\textwidth]{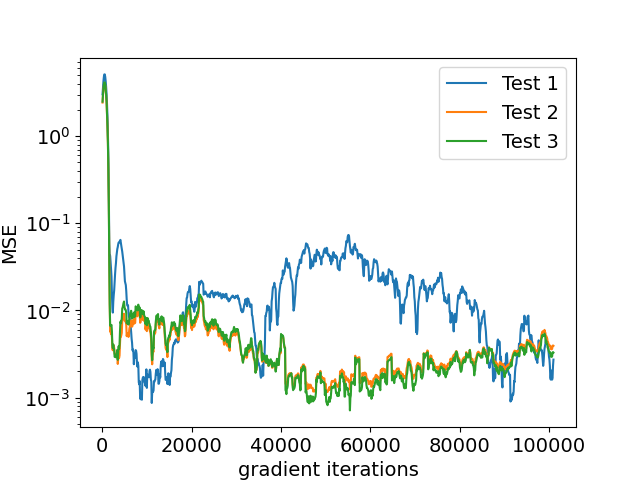}
\caption*{\tiny Test error with   cylinder  : N=50000}
\end{minipage}
    \begin{minipage}[t]{0.325\linewidth}
  \centering
 \includegraphics[width=\textwidth]{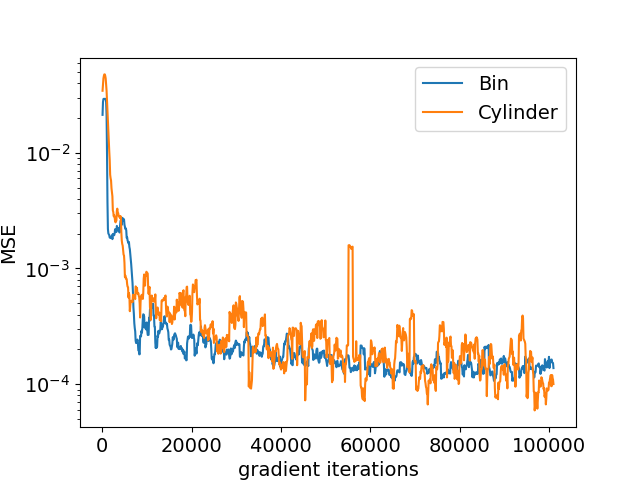}
\caption*{\tiny Training: N=250000}
\end{minipage}
\begin{minipage}[t]{0.325\linewidth}
  \centering
 \includegraphics[width=\textwidth]{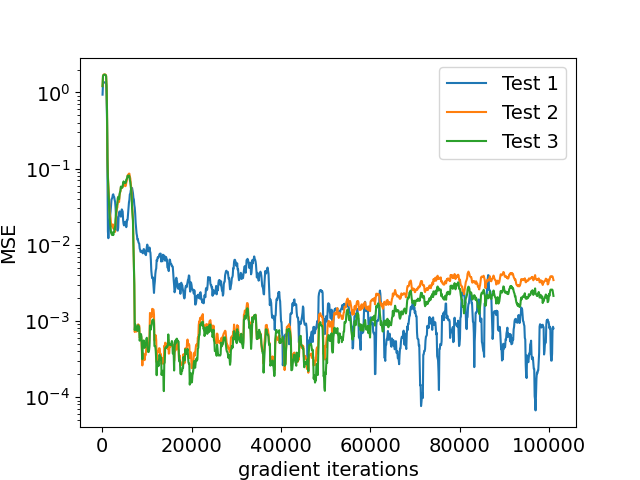}
\caption*{\tiny Test error with  bins  : N=250000}
\end{minipage}
  \begin{minipage}[t]{0.325\linewidth}
  \centering
 \includegraphics[width=\textwidth]{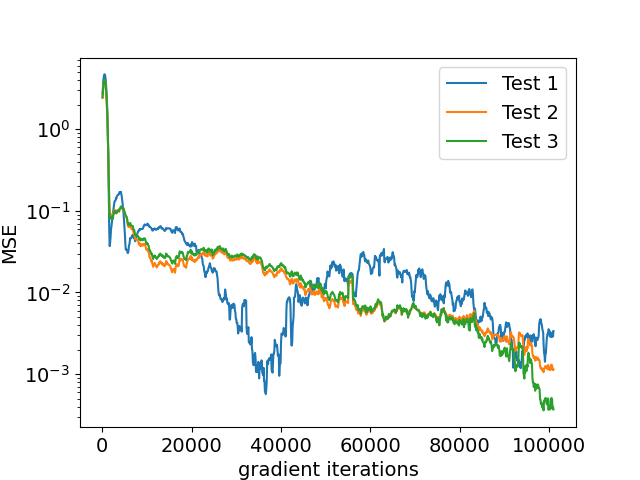}
\caption*{ \tiny Test error with   cylinder : N=250000}
\end{minipage}
     \caption{Case A: comparing bin to cylinder methods. Convergence for given distributions. 
      }
     \label{fig:CaseA}
 \end{figure}

\begin{figure}[H]
     \centering
     \begin{minipage}[t]{0.325\linewidth}
  \centering
 \includegraphics[width=\textwidth]{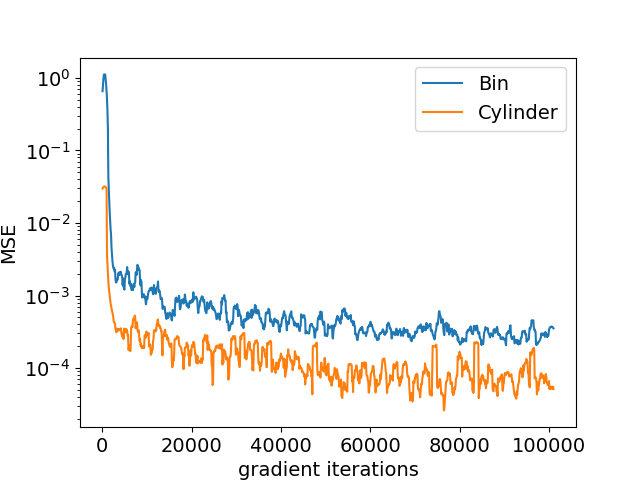}
\caption*{\tiny Training  MSE : N=10000}
\end{minipage}
\begin{minipage}[t]{0.325\linewidth}
  \centering
 \includegraphics[width=\textwidth]{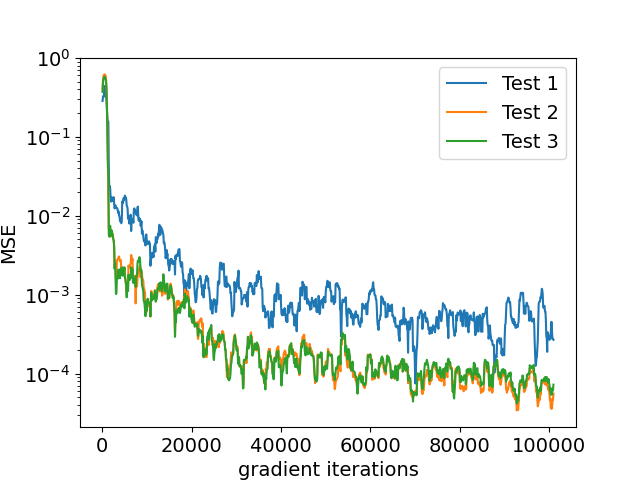}
\caption*{\tiny Test error with  bins : N=10000}
\end{minipage}
  \begin{minipage}[t]{0.325\linewidth}
  \centering
 \includegraphics[width=\textwidth]{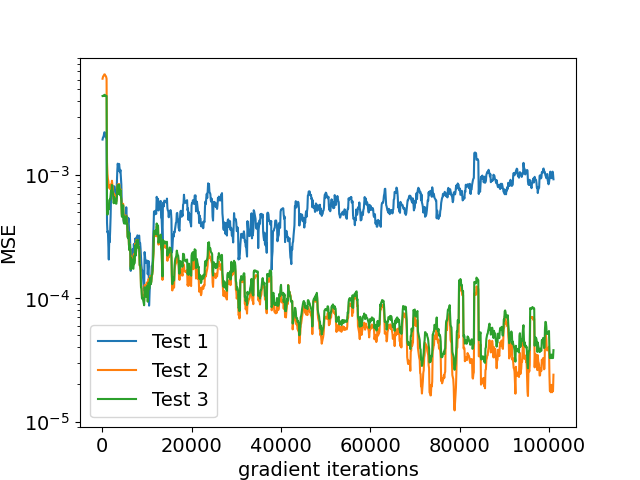}
\caption*{\tiny Test error with  cylinder  : N=10000}
\end{minipage}
     \begin{minipage}[t]{0.325\linewidth}
  \centering
 \includegraphics[width=\textwidth]{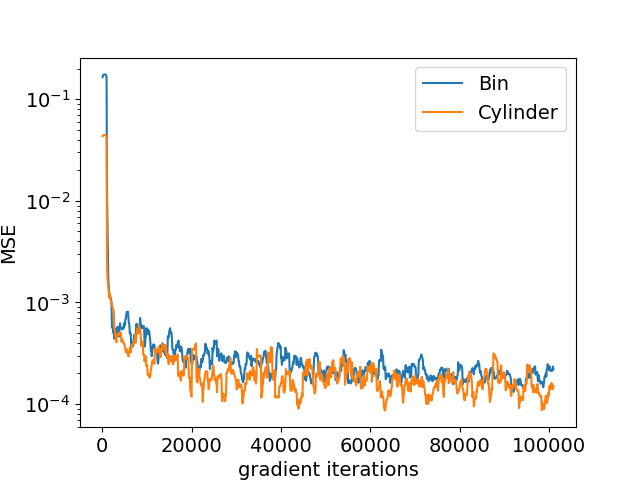}
\caption*{\tiny Training  MSE : N=50000}
\end{minipage}
\begin{minipage}[t]{0.325\linewidth}
  \centering
 \includegraphics[width=\textwidth]{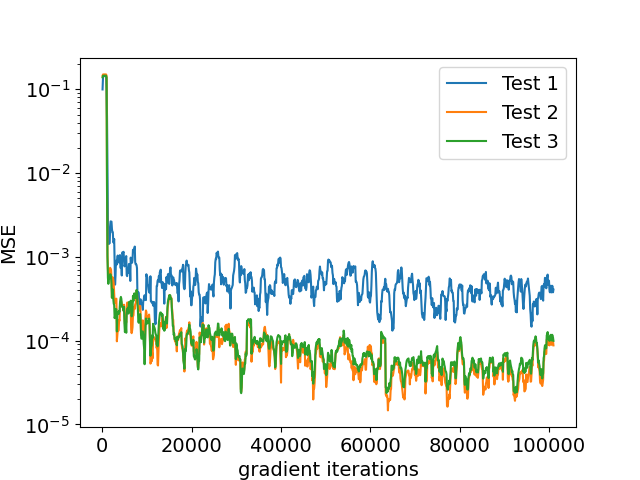}
\caption*{ \tiny Test error with  bins : N=50000}
\end{minipage}
  \begin{minipage}[t]{0.325\linewidth}
  \centering
 \includegraphics[width=\textwidth]{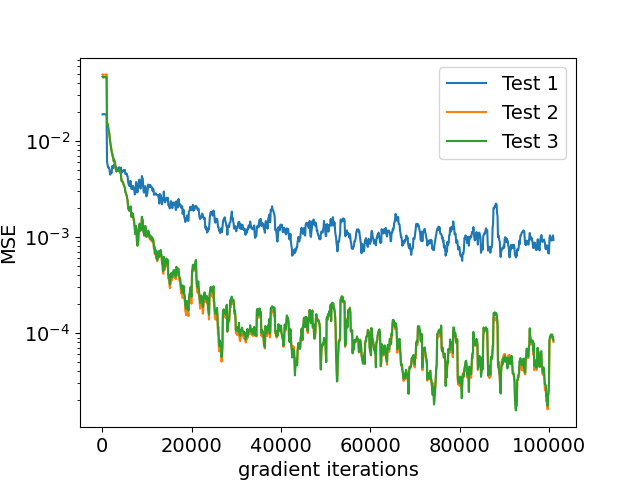}
\caption*{\tiny Test error with  cylinder  : N=50000}
\end{minipage}
    \begin{minipage}[t]{0.325\linewidth}
  \centering
 \includegraphics[width=\textwidth]{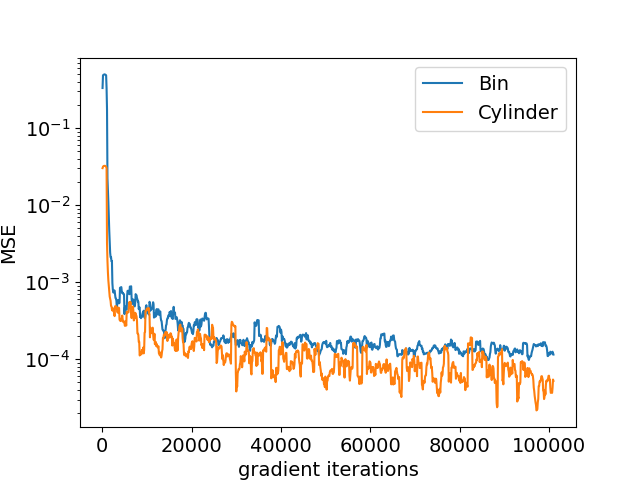}
\caption*{\tiny Training  MSE : N=250000}
\end{minipage}
\begin{minipage}[t]{0.325\linewidth}
  \centering
 \includegraphics[width=\textwidth]{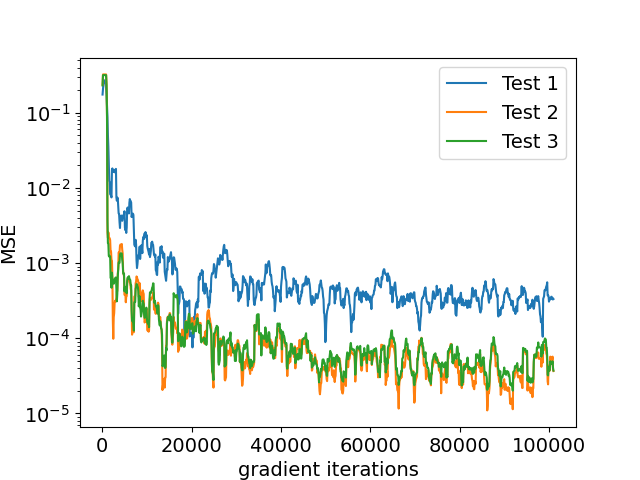}
\caption*{\tiny Test error with  bins : N=250000}
\end{minipage}
  \begin{minipage}[t]{0.325\linewidth}
  \centering
 \includegraphics[width=\textwidth]{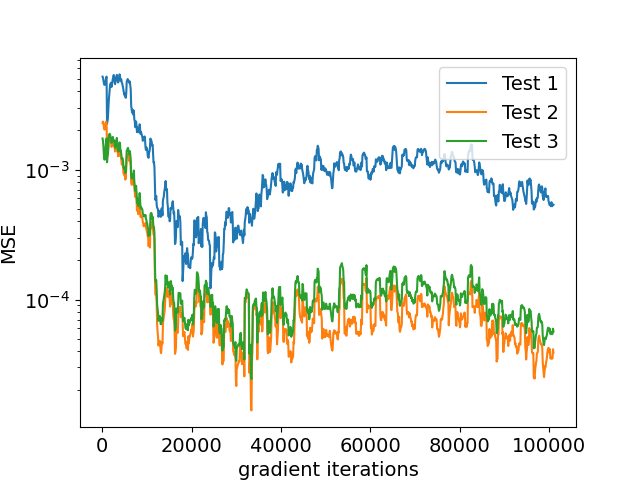}
\caption*{ \tiny Test error with  cylinder  : N=250000}
\end{minipage}
     \caption{Case B : comparing bin to cylinder methods. Convergence for given distributions.}
     \label{fig:CaseB}
 \end{figure}

On Figures \ref{fig:CaseA} and \ref{fig:CaseB}, we plot for different values of $N$,  keeping  $\hat M=20$, the training MSE obtained by the algorithms and the error associated to the  test distributions. 
Globally, the bin approximation is  more sensitive to the $N$ parameter. 
Not surprisingly, the training size  $N$ has to be taken large to get good results for both methods.
For the case B, the results on test 1 with a distribution tending to the Dirac distribution has a residual error clearly above the errors obtained with the distribution test 2 and 3.

 Finally, we consider three other mean-field functions
 \begin{itemize}
 \item [C] {\it Case C : a second-order mean-field interaction}
    \begin{align}
     V(x, \mu)  &= \;  \int \int (x - y - z )^2 \mu(\d y) \mu(\d z) \\
     \; &=  \; x^2  - 4 x\bar\mu + 2 \E_{X\sim\mu}[X^2] +  2 |\bar\mu|^2. 
    \end{align} 
     \item[D.]  {\it Case D: median function}
     \begin{align}
    V(x, \mu)  & = \;  \int | x- y| \mu(dy) \; = \;  \E_{X \sim \mu}|x- X|. 
    \end{align}
    \item[E.] {\it  Case E : cumulative distribution  function} 
      \begin{align}
    V(x, \mu)  & = \;  \mu(-\infty,x] \; = \;   \E_{X \sim \mu}[ 1_{X \le x}] 
    \end{align}
 \end{itemize}

\begin{figure}[H]
     \centering
    \begin{minipage}[t]{0.325\linewidth}
  \centering
 \includegraphics[width=\textwidth]{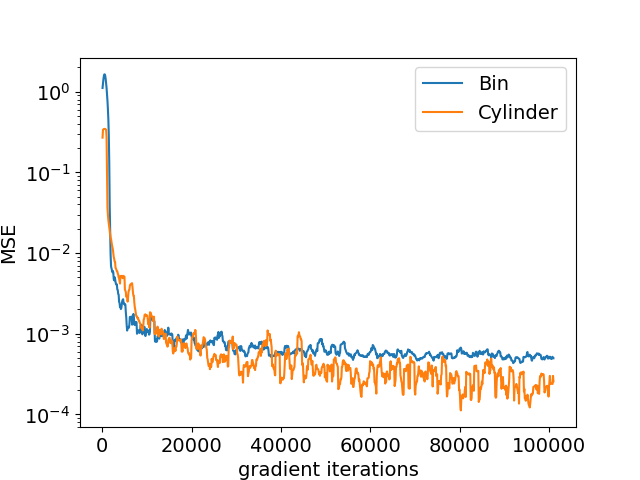}
\caption*{ Training  MSE }
\end{minipage}
\begin{minipage}[t]{0.325\linewidth}
  \centering
 \includegraphics[width=\textwidth]{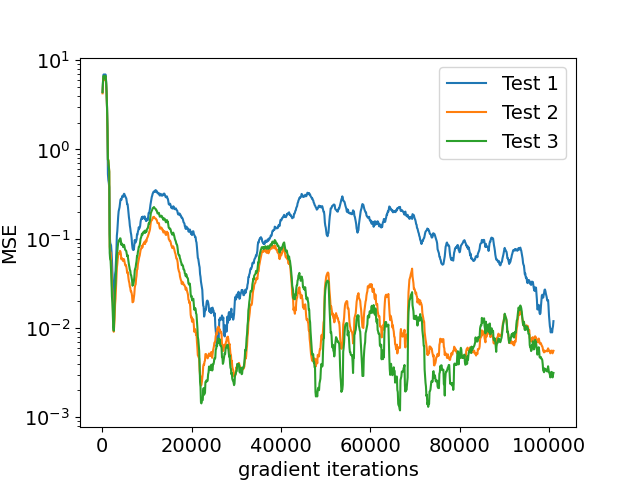}
\caption*{Test error with  bins }
\end{minipage}
  \begin{minipage}[t]{0.325\linewidth}
  \centering
 \includegraphics[width=\textwidth]{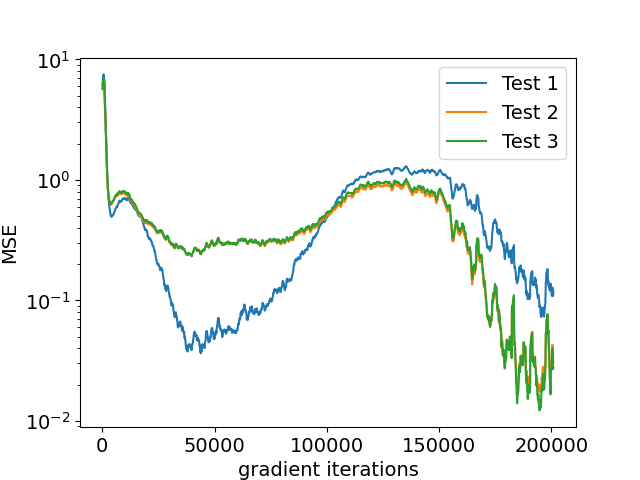}
\caption*{ \red{Test error with}  cylinder }
\end{minipage}
     \caption{Case C : comparing bin to cylinder networks with $N=250000$.}
     \label{fig:CaseC}
 \end{figure}

On Figure \ref{fig:CaseC}, we keep the same hyperparameters as  for case A and B.
For the more complex regular function of case C, the convergence is clearly more difficult to obtained on the test distributions even if the training MSE converges accurately. Results obtained by the bin network are better than the results obtained using the cylinder network. Increasing the number of layers or neurons  does not improve the results.

\begin{figure}[H]
     \centering
    \begin{minipage}[t]{0.325\linewidth}
  \centering
 \includegraphics[width=\textwidth]{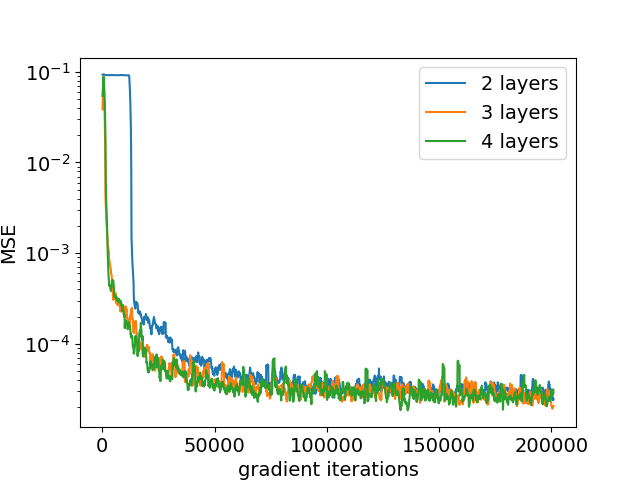}
\caption*{\red{Training}  MSE }
\end{minipage}
\begin{minipage}[t]{0.325\linewidth}
  \centering
 \includegraphics[width=\textwidth]{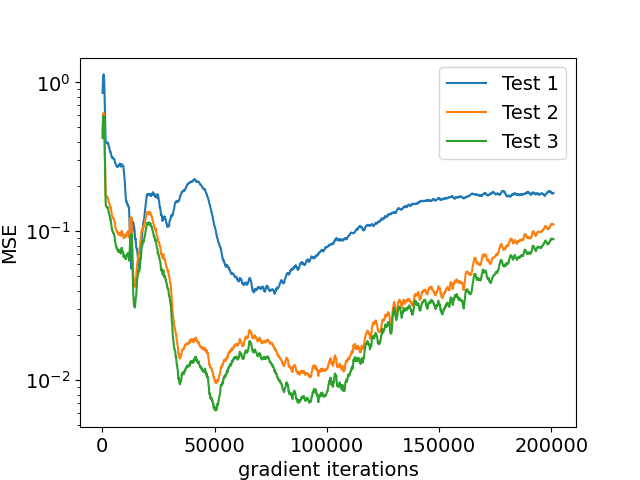}
\caption*{ \red{Test error with}  2 layers}
\end{minipage}
  \begin{minipage}[t]{0.325\linewidth}
  \centering
 \includegraphics[width=\textwidth]{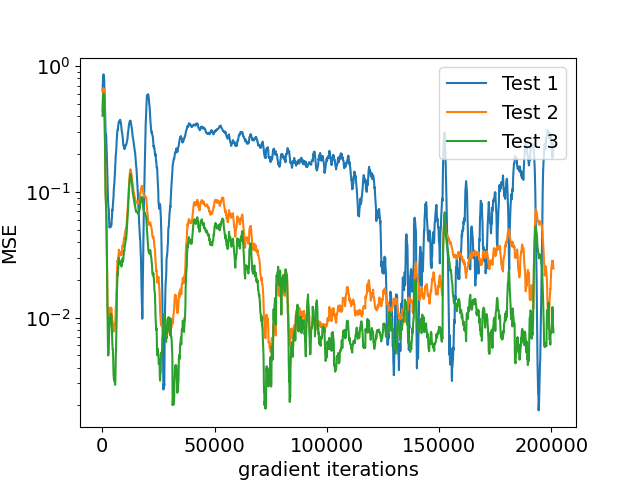}
\caption*{\red{Test error with}  4 layers }
\end{minipage}
     \caption{Case  D for bins : results depending on the number of layers with $N=250000$, $20$ neurons.}
     \label{fig:CaseDBins}
 \end{figure}

\begin{figure}[H]
     \centering
    \begin{minipage}[t]{0.325\linewidth}
  \centering
 \includegraphics[width=\textwidth]{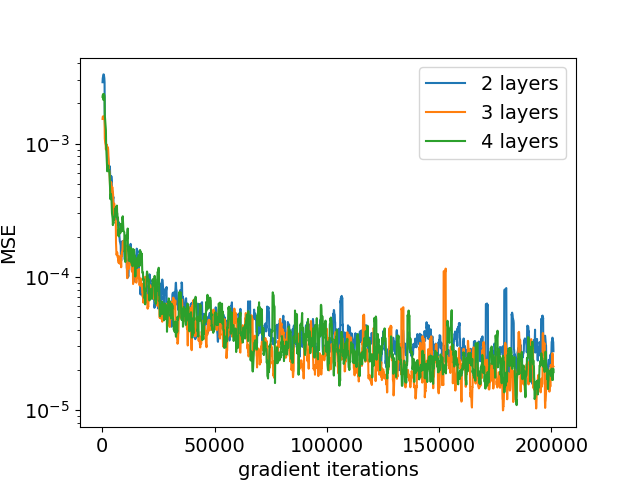}
\caption*{\small \red{Training}  MSE }
\end{minipage}
\begin{minipage}[t]{0.325\linewidth}
  \centering
 \includegraphics[width=\textwidth]{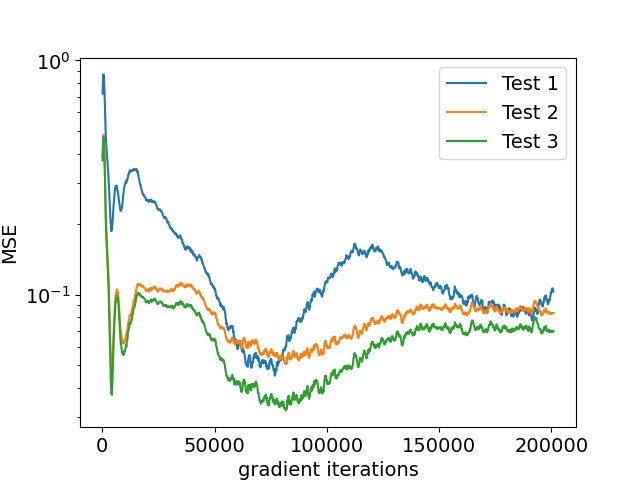}
\caption*{\small \red{Test error with}  2 layers}
\end{minipage}
  \begin{minipage}[t]{0.325\linewidth}
  \centering
 \includegraphics[width=\textwidth]{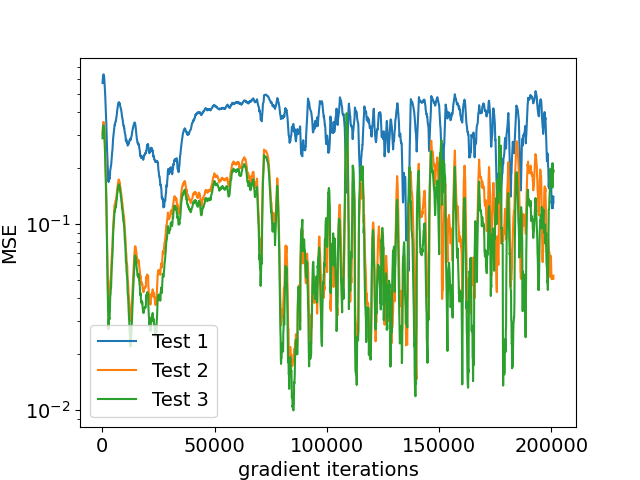}
\caption*{ \small \red{Test error with}  4 layers }
\end{minipage}
     \caption{Case  D for cylinder : results depending on the number of layers with  $N=250000$, $20$ neurons for the two networks.}
     \label{fig:CaseDCylinder}
 \end{figure}

\red{ On Figures \ref{fig:CaseDBins} and \ref{fig:CaseDCylinder}, we give the results obtained for case D by the two methods using $20$ neurons and letting the number of layers increase. The global training MSE does not seem to improve as we increase the number of layers. Again the first test function gives the worst results. A high accuracy is hard to obtained on the test distributions even if a high accuracy of the global MSE is reached.
 As we increase the number of layers, the results on the test distributions are more erratic.}


\begin{figure}[H]
     \centering
    \begin{minipage}[t]{0.325\linewidth}
  \centering
 \includegraphics[width=\textwidth]{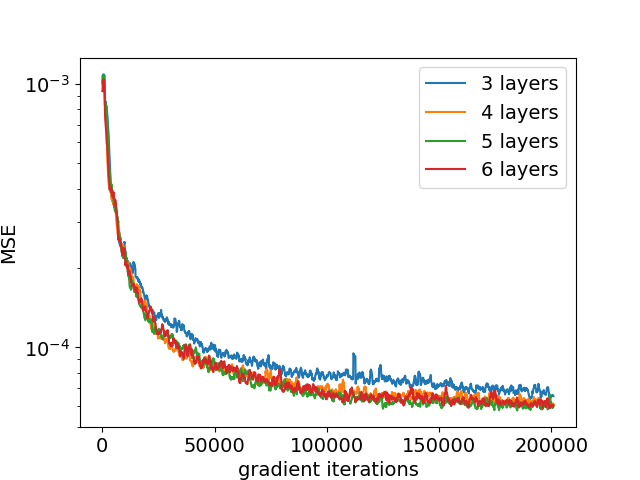}
\caption*{\small \red{Training}  MSE }
\end{minipage}
\begin{minipage}[t]{0.325\linewidth}
  \centering
 \includegraphics[width=\textwidth]{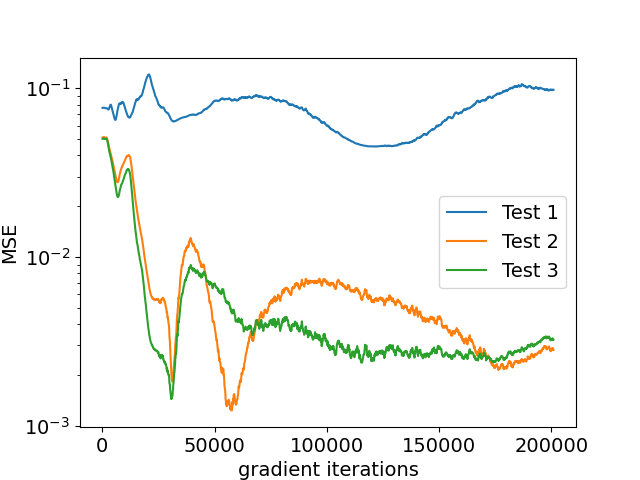}
\caption*{\small \red{Testing}   2 layers }
\end{minipage}
  \begin{minipage}[t]{0.325\linewidth}
  \centering
 \includegraphics[width=\textwidth]{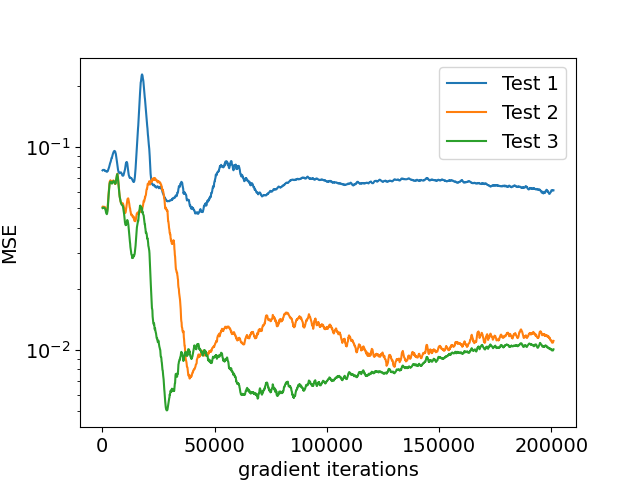}
\caption*{ \small \red{Testing}  4 layers}
\end{minipage}
     \caption{Case E for bins : results depending on the number of layers $N=250000$, $20$ neurons.}
     \label{fig:CaseEBins}
 \end{figure}

\begin{figure}[H]
     \centering
    \begin{minipage}[t]{0.325\linewidth}
  \centering
 \includegraphics[width=\textwidth]{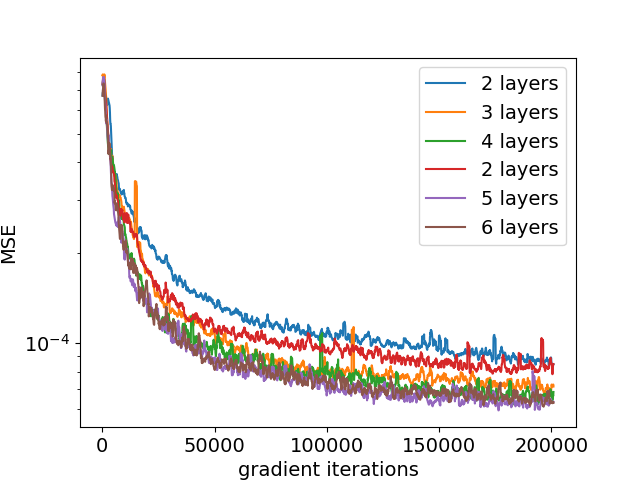}
\caption*{\small \red{Training}  MSE }
\end{minipage}
\begin{minipage}[t]{0.325\linewidth}
  \centering
 \includegraphics[width=\textwidth]{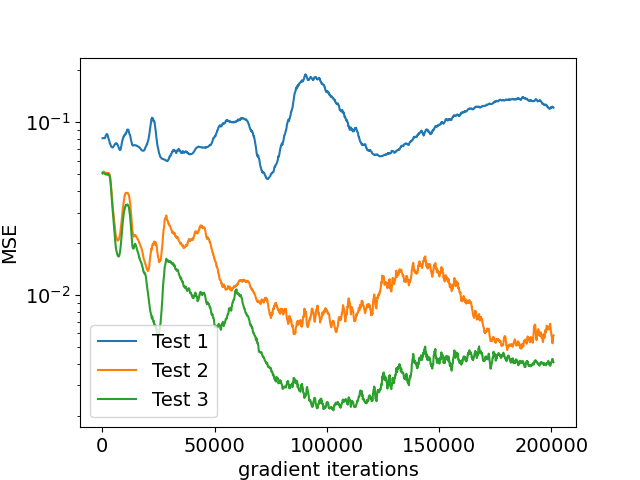}
\caption*{\small \red{Test error with}  2 layers}
\end{minipage}
  \begin{minipage}[t]{0.325\linewidth}
  \centering
 \includegraphics[width=\textwidth]{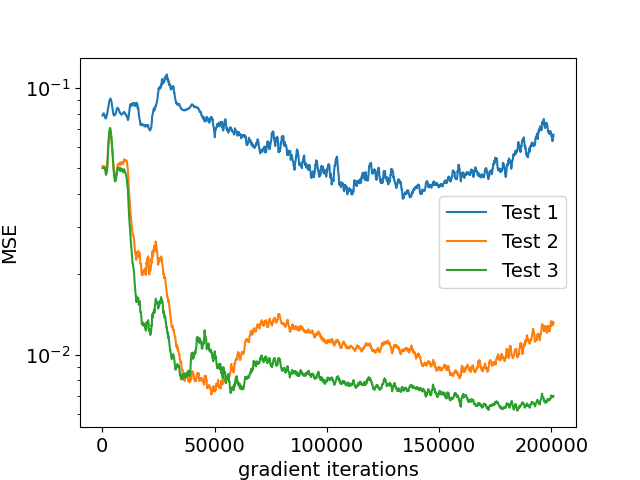}
\caption*{ \small \red{Test error with}  4 layers }
\end{minipage}
     \caption{Case  E for cylinder : results depending on the number of layers with  $N=250000$, $20$ neurons for the two networks.}
     \label{fig:CaseECylinder}
 \end{figure}

 \red{
On Figures \ref{fig:CaseEBins} and \ref{fig:CaseECylinder}, we test the influence of the number of layers of the networks on case E using $20$ neurons. For this very irregular function, increasing the number of layers improve the training MSE and the results on the test distributions but only marginally. Again results obtained on the first test distribution are not as good as the results obtained for the two other distributions.
}

\section{Algorithms for dynamic mean-field problems} \label{sec:dynMF}

In dynamic mean-field problems (over finite horizon), like mean-field control/game, the solution (value function, control) is time-dependent,  and function of some state process and its probability distribution.  It is then defined on $\Tc\times\R^d\times\Pc_2(\R^d)$, where $\Tc$ is 
an interval of the form $[0,T]$ in a continuous-time problem, or a discrete time grid $\Tc$ $=$ $\{t_i, i=0,\ldots,N_T\}$ in a discrete-time problem. Denoting by $U$ this time-dependent mean-field function, we aim to approximate by learning the time-dependent functional 
\begin{align} \label{defVcMF} 
\Uc(t) : \mu \in \Pc_2(\R^d) & \longmapsto \;  U(t,\cdot,\mu) \in L^2(\mu),  \quad \mbox{ for } t \in \Tc. 
\end{align} 
The solution $U$ is typically characterized via a Master Bellman PDE, a dynamic programming formula,  or a 
Backward stochastic differential equations of McKean-Vlasov type, and after time discretization (in the case of a continuous-time problem) on a grid $\{t_i,i=0,\ldots,N_T\}$,  one can design algorithms for learning the operators $\Uc_i$ $:=$ $\Uc(t_i)$, $i$ $=$ $0,\ldots,N_T$. These machine learning algorithms are either of global or local type, and are described in the next paragraphs.

\subsection{Local algorithms}

In backward recursion  local algorithm, arising from dynamic programming,  
given an approximation at time $t_{i+1}$ of  the mean-field operator $\Uc_{i+1}$ (e.g. the value function and or the feedback control)  by a mean-field neural network function $\widehat\Nc_{i+1}$ 
as described in the previous section, we aim to learn at time $t_i$  a mean-field neural network function $\Nc_\theta$ by minimizing over $\theta$ a loss function in the form
\beqs
L_i(\theta) &=& \E \big[   H \big( X_i,\mu_i,\Nc_\theta(\mu_i)(X_i),\widehat \Nc_{i+1}(\mu_{i+1})(X_{i+1}) \big) \big],   
\enqs
for some function $H$, and we then update $\widehat\Nc_{i}$ $=$ $\Nc_{\hat\theta_i}$ where $\hat\theta_i$ is the resulting optimal parameter from this  minimization problem. 
In the above expectation for applying SGD, $\mu_i$ is sampled according to the data generation as described in the previous section, $X_i$ is sampled according to $\mu_i$,  $X_{i+1}$ is given by a  dynamics in the form: 
\beqs
X_{i+1} &=&  F_i(X_i,\mu_i,\Nc_\theta(\mu_i)(X_i),\eps_{i+1}),
\enqs
and $\mu_{i+1}$ is the law of $X_{i+1}$.  In practice, $\mu_{i+1}$ has to be estimated/approximated from samples of $X_{i+1}$, and the suitable method will be chosen depending on the adopted  class of mean-field neural networks. 
\begin{enumerate}
\item {\it Bin density-based neural network}:  $\Nc_\theta(\mu)$ $=$ $\Phi_\theta(.,\bop^\mu)$. We sample probability measures $\mu_i^{(m)}$ $=$ $\Lc_D(\bop^{(m)})$ in $\Dc_2(\R)$ 
from samples \\ $\bop^{(m)}$ $=$ $(p_k^{(m)})_{k\in\llbracket 1,K\rrbracket}$, $m$ $=$ $1,\ldots,M$, in $\Dc_K$, and then approximate the loss function $L_i$ as
\begin{footnotesize}
\beqs
 & & L_i(\theta) \\
 & \simeq &  \frac{1}{MN} \sum_{m=1}^M \sum_{n=1}^{N}    H \big(X_i^{(m),(n)},\mu_i^{(m)},\Phi_\theta(X_i^{(m),(n)},\bop^{(m)}),\Phi_{\hat\theta_{i+1}}(X_{i+1}^{(m),(n)},\hat\bop^{(m)}) \big) 
\enqs
\end{footnotesize}
where $X_i^{(m),(n)}$, $n$ $=$ $1,\ldots,N$, are sampled from $\Lc_D(\bop^{(m)})$,  
\beqs
X_{i+1}^{(m),(n)} &=& F(X_i^{(m),(n)},\Lc_D(\bop^{(m)}),\Phi_\theta(X_i^{(m),(n)},\bop^{(m)}),\eps_{i+1}),
\enqs
and $\hat\bop^{(m)}$ $=$ $(\hat p_k^{(m)})_{k \in \llbracket 1,K\rrbracket}$ are the estimated density bins in $\Dc_K$ of  $X_{i+1}^{(m),(n)}$ (truncated on $\Kc$ $=$ $[x_0,x_K]$), namely: 
\beqs
\hat p_{k}^{(m)} &=& \frac{ \#  \{ n \in \llbracket 1,N\rrbracket:  {\rm Proj}_\Kc(X_{i+1}^{(m),(n)}) \in \mbox{Bin}(k) \} }{Nh}, \quad k=1,\ldots,K,
\enqs
where ${\rm Proj}_{\Kc}(.)$ is the projection on $\Kc$. 
\item {\it Cylindrical neural network}:
$\Nc_\theta(\mu)$ $=$ $\Psi_\theta(.,<\varphi_\theta,\mu>)$.  We sample probability measures  $\mu^{(m)}$ $=$ $\Lc_D(\bop^{(m)})$, $m=1,\ldots,M$,
 so that the loss function is approximated as
  \begin{footnotesize}
  \beqs
  & & L_i(\theta) \\
  & \simeq &  \frac{1}{MN} \sum_{m=1}^M \sum_{n=1}^{N}  H \big(X_i^{(m),(n)},\mu_i^{(m)}, \Psi_\theta(X_i^{(m),(n)},  \hat \E[\varphi_\theta(X_i^{(m)})]), 
  \Psi_{\hat\theta_{i+1}}( X_{i+1}^{(m),(n)},  \\
  && \quad \quad \quad  \quad \quad \quad  \hat\E[ \varphi_{\hat\theta_{i+1}}(  X_{i+1}^{(m)} ) ] ) \big), 
 \enqs
 \end{footnotesize}
 where $X_i^{(m),(n)}$, $n$ $=$ $1,\ldots,N$, are sampled from $X_i^{(m)}$ $\sim$ $\Lc_D(\bop^{(m)})$,  $X_{i+1}^{(m),(n)}$, $n$ $=$ $1,\ldots,N$, are sampled as
\beqs
X_{i+1}^{(m),(n)} &=& F_i(X_i^{(m),(n)},\mu_i^{(m)}, \Psi_\theta(X_i^{(m),(n)},  \hat \E[\varphi_\theta(X_i^{(m)})]),\eps_{i+1}^{(n)}),
\enqs
with $\hat\E[.]$ denoting the empirical expectation: 
\beqs
 \hat \E[\varphi_\theta(X_i^{(m)})] &=&  \frac{1}{N} \sum_{n=1}^N \varphi_\theta(X_i^{(m),(n)}), \\ \hat\E[ \varphi_{\hat\theta_{i+1}}(  X_{i+1}^{(m)} ) ] \; & = & \;  \frac{1}{N} \sum_{n=1}^N  \varphi_{\hat\theta_{i+1}}(X_{i+1}^{(m),(n)}). 
\enqs
\end{enumerate}

\subsection{Global algorithms}

In global algorithms, we approximate at any time $t_i$, $i$ $=$ $0,\ldots,N_T$, the mean-field operators $\Uc_i$ by mean-field networks $\Nc_{\theta_i}$ that are learned simultaneously by minimizing over $\botheta$ $=$ $(\theta_i)_i$  a global loss function in the form 
\beqs
L(\botheta) &=& \E \Big[ \sum_{i=0}^{N_T-1} \ell_i\big(X_i,\mu_i,\Nc_{\theta_i}(\mu_i)(X_i) \big) + g(X_{N_T},\mu_{N_T}) \Big]
\enqs  
for some loss functions $\ell_i$, and $g$.  In the above expectation for applying SGD,  $\mu_0$ is sampled according to the data generation as described in the previous section, $X_0$ is sampled according to $\mu_0$, and for $i$ $=$ $0,\ldots,N_T-1$, $X_{i+1}$ are given by a dynamics in the form 
\beqs
X_{i+1} &=&  F_i(X_i,\mu_i,\Nc_{\theta_i}(\mu_i)(X_i),\eps_{i+1}),
\enqs
where $\mu_i$ is the law of $X_i$.   In practice, for $i$ $=$ $1,\ldots,N_T$, $\mu_{i}$ has to be estimated/approximated from samples of $X_{i}$, and the suitable method will be chosen depending on the adopted  class of mean-field neural networks. 
\begin{enumerate}
\item  {\it Bin density-based neural network}:  $\Nc_\theta(\mu)$ $=$ $\Phi_\theta(.,\bop^\mu)$. We sample probability measures $\mu_0^{(m)}$ $=$ $\Lc_D(\bop^{(m)})$ in $\Dc_2(\R)$ 
from samples $$\bop^{(m)} = (p_k^{(m)})_{k\in\llbracket 1,K\rrbracket}, m = 1,\ldots,M,$$ in $\Dc_K$, and then approximate the  global loss function as
\beqs
L(\botheta) & \simeq & \frac{1}{MN} \sum_{m=1}^M \sum_{n=1}^N  \Big[ \ell_0(X_0^{(m),(n)},\mu_0^{(m)},\Phi_{\theta_0}(X_0^{(m),(n)},\bop^{(m)})  \big) \\
& & \quad + \; \sum_{i=1}^{N_T-1}  \ell_i\big(X_i^{(m),(n)},\hat\mu_i^{(m)},\Phi_{\theta_i}(X_i^{(m),(n)}, \hat\bop_i^{(m)}  \big)  \\
& &  \quad + g(X_{N_T}^{(m),(n)},\hat\mu_{N_T}^{(m)}) \Big] 
\enqs
where $X_0^{(m),(n)}$, $n$ $=$ $1,\ldots,N$ are sampled from $X_0^{(m)}$ $\sim$ $\mu_0^{(m)}$, for $i$ $=$ $0,\ldots,N_T-1$,  $X_{i+1}^{(m),(n)}$, $n$ $=$ $1,\ldots,N$ are sampled as 
\beqs
X_{i+1}^{(m),(n)} &=& F_i(X_i^{(m),(n)},\hat\mu_i^{(m)},\Phi_{\theta_i}(X_i^{(m),(n)},\hat\bop_i^{(m)}),\eps_{i+1}^{(n)}),   
\enqs
with $\hat\mu_i^{(m)}$ $=$ $\Lc_D(\hat\bop_i^{(m)})$, $\hat\bop_0^{(m)}$ $=$ $\bop^{(m)}$, and $\hat\bop_i^{(m)}$ $=$ $(\hat p_{i,j}^{(m)})_{j \in \llbracket 1,K\rrbracket}$ are the estimated density weights in $\Dc_K$ of  $X_{i}^{(m),(n)}$, $i$ $=$ $1,\ldots,N_T$ (truncated on $\Kc$ $=$ $[x_0,x_K]$), namely: 
\beqs
\hat p_{i,j}^{(m)} &=& \frac{ \#  \{ n \in \llbracket 1,N\rrbracket:  {\rm Proj}_\Kc(X_{i}^{(m),(n)}) \in \mbox{Bin}(j) \} }{Nh_j}, \quad j=1,\ldots,K,
\enqs
where ${\rm Proj}_{\Kc}(.)$ is the projection on $\Kc$. 
\item {\it Cylindrical neural network}: $\Nc_\theta(\mu)$ $=$ $\Psi_\theta(.,<\varphi_\theta,\mu>)$.
We sample probability measures $\mu_0^{(m)}$ say according to Bin density measures $\Lc_D(\bop^{(m)})$, and then minimize over the parameters $\botheta$ $=$ $(\theta_i)$ the approximate global loss function 
\beqs
L(\botheta) & \simeq & \frac{1}{MN} \sum_{m=1}^M \sum_{n=1}^N  \Big[ \ell_0(X_0^{(m),(n)},\mu_0^{(m)},\Psi_{\theta_0}(X_0^{(m),(n)}, \hat \E[\varphi_\theta(X_0^{(m)})])  \big) \\
& & \quad + \; \sum_{i=1}^{N_T-1}  \ell_i\big(X_i^{(m),(n)},\hat\mu_i^{(m)},\Psi_{\theta_i}(X_i^{(m),(n)},  \hat \E[\varphi_\theta(X_i^{(m)})] \big)  \\
& & \quad + g(X_{N_T}^{(m),(n)},\hat\mu_{N_T}^{(m)}) \Big], 
\enqs
where $X_0^{(m),(n)}$, $n$ $=$ $1,\ldots,N$ are sampled from $X_0^{(m)}$ $\sim$ $\mu_0^{(m)}$, for $i$ $=$ $0,\ldots,N_T-1$,  $X_{i+1}^{(m),(n)}$, $n$ $=$ $1,\ldots,N$ are sampled as
\beqs
X_{i+1}^{(m),(n)} &=& F_i(X_i^{(m),(n)},\hat\mu_i^{(m)},\Psi_{\theta_i}(X_i^{(m),(n)}, \hat \E[\varphi_\theta(X_i^{(m)})]),\eps_{i+1}),   
\enqs
$\hat\mu_i^{(m)}$ $=$ $\Lc_D(\hat\bop_i^{(m)})$, $\hat\bop_0^{(m)}$ $=$ $\bop^{(m)}$, and $\hat\bop_i^{(m)}$ $=$ $(\hat p_{i,j}^{(m)})_{j \in \llbracket 1,K\rrbracket}$ are the estimated density weights in $\Dc_K$ of  $X_{i}^{(m),(n)}$, $i$ $=$ $1,\ldots,N_T$ (truncated on $\Kc$ $=$ $[x_0,x_K]$), namely: 
\beqs
\hat p_{i,j}^{(m)} &=& \frac{ \#  \{ n \in \llbracket 1,N\rrbracket:  {\rm Proj}_\Kc(X_{i}^{(m),(n)}) \in \mbox{Bin}(j) \} }{Nh_j}, \quad j=1,\ldots,K,
\enqs
with $\hat\E[.]$ denoting the empirical expectation: 
\beqs
 \hat \E[\varphi_\theta(X_i^{(m)})] &=&  \frac{1}{N} \sum_{n=1}^N \varphi_\theta(X_i^{(m),(n)}). 
\enqs
\end{enumerate}

\begin{Remark}
For global algorithms, we can avoid the approximation of the mean-field function at each single date $t_i$, $i$ $=$ $0,\ldots,N_T-1$, by learning directly the mean field function  which also takes the  time as argument. 
Hence, instead of having a different mean-field neural network $\Nc_{\theta_i}$ for each date $t_i$, we learn with  a single time dependent neural network $\Nc(t,.)$, which is  used for all dates, as illustrated on an example in the next section.  
This gives generally more stable results, see e.g.  \cite{chan2019machine}.
\end{Remark}

\subsection{A toy example of semi-linear PDE on Wasserstein space}

Let us consider the linear differential operator on $[0,T]\times\R\times\Pc_2(\R)$ associated to the mean-field stochastic differential equation (SDE): 
\beqs
\d X_t &=& \kappa( \E[X_t] - X_t) \d t + \sigma  \d W_t, 
\enqs 
for some positive constants $\kappa$, $\sigma$, and given by
\beqs
\Lc  v(t,x,\mu) &=& \Dt{v}(t,x,\mu) +  \kappa(\bar\mu-x) \Dx{v}(t,x,\mu) + \frac{1}{2} \sigma^2 \Dxx{v}(t,x,\mu) \\
& & \quad \quad + \;  \E_{\xi\sim\mu} \big[ \kappa(\bar\mu-\xi) \partial_\mu v(t,x,\mu)(\xi) +  \frac{1}{2} \sigma^2 \partial_{x'}\partial_\mu v(t,x,\mu)(\xi) \big], 
\enqs  
where $x'$ $\mapsto$ $\partial_\mu v(t,x,\mu)(x')$ is the Lions derivative of $\mu$ $\mapsto$ $v(t,x,\mu)$ (see \cite{cardel19}).  

Given a $C^2$ function $w$ on $\R^d$ with quadratic growth condition, let us define the function $f$ on $[0,T]\times\R\times\Pc_2(\R)\times\R$ by 
\begin{align} 
f(t,x,\mu,y) &= \;  e^{T-t}   \E_{\xi\sim\mu} \Big[ (w -  \sigma^2 D_{xx}w)(x-\xi) +  \kappa(x-\xi) D_x w(x-\xi) \Big] \\
& \quad  \quad  - \; a  \Big( \E_{\xi\sim\mu}\big[ e^{T-t} w(x-\xi) \big] \Big)^2 + a y^2, \label{deff} 
\end{align} 
for some positive constant $a$, and  consider the semi-linear PDE on $[0,T]\times\R\times\Pc_2(\R)$: 
\begin{equation} \label{eq:semiPDE} 
\left\{
\begin{array}{rcl}
\Lc v + f(t,x,\mu,v) &=& 0,  \quad \quad \quad (t,x,\mu) \in [0,T)\times\R\times\Pc_2(\R), \\
v(T,x,\mu) &=&  g(x,\mu),  \quad (x,\mu) \in \R\times\Pc_2(\R), 
\end{array}
\right.
\end{equation} 
where $g(x,\mu)$ $:=$ $\E_{\xi\sim\mu}[w(x-\xi)]$.  By construction, the solution to the PDE \eqref{eq:semiPDE} is explicitly given by $v(t,x,\mu)$ $=$ $e^{T-t} \E_{\xi\sim\mu}[w(x-\xi)]$, and this will serve as benchmark for evaluating the accuracy of our 
algorithms in the numerical resolution of the PDE  \eqref{eq:semiPDE}. 

Recall that the PDE  \eqref{eq:semiPDE} has the following  probabilistic representation: by considering the pair of processes  $(Y,Z)$ given by 
\beqs
Y_t \;  = \; v(t,X_t,\P_{X_t}), & &  Z_t \; = \;  \sigma \Dx v(t,X_t,\P_{X_t}), 
\quad 0 \leq t\leq T, 
\enqs
we see  by It\^o's formula that it satisfies  the Backward SDE
\begin{align} \label{BSDEsemi} 
\d Y_t &= \:   - f(t,X_t,\P_{X_t},Y_t) \d t  + Z_t \d W_t, \quad Y_T =  
g(X_T,\P_{X_T}).  
\end{align}

\vspace{2mm}

\paragraph{Local Algorithms.}  We consider a time grid $\Tc$ $=$ $\{t_i,i=0,\ldots,N_T\}$ of $[0,T]$ with mesh size $\Delta t_i$ $=$ $t_{i+1}-t_i$, and consider two local algorithms for approximating $v$ on $\Tc\times\R\times\Pc_2(\R)$. 
In the first approach, we start from the expectation representation arising from \eqref{BSDEsemi}:  
\begin{footnotesize}
\beqs
v(t_i,X_{t_i},\P_{X_{t_i}}) &=& \E \Big[ v(t_{i+1},X_{t_{i+1}},\P_{X_{t_{i+1}}}) +  \int_{t_i}^{t_{i+1}} f(s,X_s,\P_{X_s},v(s,X_s,\P_{X_s}) ) \d s \big| X_{t_i} \Big], 
\enqs
\end{footnotesize}
which leads to the backward regression algorithm: starting from $\hat\Uc_{N_T}(\mu)(x)$ $=$ $g(x,\mu)$, we approximate $v$ at any time $t_i$, $i$ $=$ $N_T-1,\ldots,0$, by mean-field neural networks $\Uc_{\theta_i}$, and minimize the local loss regression function
\begin{footnotesize}
\beqs
L_i^{R}(\theta_i) &=&  \E \Big| \hat\Uc_{i+1}(\mu_{i+1})(X_{i+1}) - \Uc_{\theta_i}(\mu_i)(X_i) + f(t_i,X_i,\mu_i,\Uc_{\theta_i}(\mu_i)(X_i)) \Delta t_i  \Big|^2, 
\enqs 
\end{footnotesize}
by updating $\hat\Uc_i$ $=$ $\Uc_{\hat\theta_i}$ with $\hat\theta_i$ the resulting ``optimal" parameter, and where we sample  $\mu_i$, $X_i$ $\sim$ $\mu_i$, with $(X_i)_i$ given by the Euler scheme of the mean-field SDE: 
\beqs
X_{i+1} &=& X_i + \kappa(\bar\mu_i - X_i) \Delta t_i + \sigma \Delta W_{t_i}, \quad \Delta W_{t_i} := W_{t_{i+1}} - W_{t_i}. 
\enqs

Alternately, by relying directly on the time discretization of the BSDE \eqref{BSDEsemi}, and follo\-wing the idea in \cite{hure2020deep},  we approximate $v$ and its gradient $\sigma D_x v$ at any time $t_i$ by mean-field neural networks $\Uc_{\theta_i}$ and $\Zc_{\theta_i}$,  and minimize the loss function
\beqs
L_i^{BSDE}(\theta_i) &=& \E \Big| \hat\Uc_{i+1}(\mu_{i+1})(X_{i+1}) - \Uc_{\theta_i}(\mu_i)(X_i) \\
& & \quad \quad + \;  f(t_i,X_i,\mu_i,\Uc_{\theta_i}(\mu_i)(X_i)) \Delta t_i - \Zc_{\theta_i}(\mu_i)(X_i) \Delta W_{t_i} \Big|^2. 
\enqs

\paragraph{Global Algorithms.} 
We propose two global methods. In the first one, we approximate $v$ at any time $t_i$, $i$ $=$ $0,\ldots,N_T$, by  a time dependent mean-field neural networks $\Uc_{\theta}(t,.)$, and minimize over $\theta$  the global loss regression function: 
\begin{footnotesize}
\beqs
L^R(\theta) &=& 
\E \Big[ \big|  g(X_{N_T},\mu_{N_T})  - \Uc_{\theta}( t_{N_T},\mu_{N_T})(X_T)\big|^2  \\
& & + \;  \sum_{i=0}^{N_T-1} \big| \Uc_{\theta}(t_{i+1},\mu_{i+1})(X_{i+1}) - \Uc_{\theta}(t_i,\mu_i)(X_i)  
\\  & & + f(t_i,X_i,\mu_i,\Uc_{\theta}(t_i,\mu_i)(X_i)) \Delta t_i  \big|^2 \Big].  
\enqs
\end{footnotesize}
Alternately, following the idea in \cite{weinan2017deep}, we approximate $v$ at time $t$ $=$ $0$  by a  mean-field neural networks $\Uc_{\bar \theta}(.)$, and its gradient  $\sigma D_x v$ at any time $t_i$ by a time dependent mean-field neural networks $\Zc_{\tilde \theta}(t,.)$ 
by minimizing over $\theta = (\bar \theta, \tilde \theta) $  the global loss function:  
\beqs
L^{BSDE}(\theta) &=& \E \Big| Y_{N_T}^\theta - g(X_T,\P_{X_T}) \Big|^2,\red
\enqs
where $Y^\theta$ is  given by 
\begin{footnotesize}
\beqs
Y_{i+1}^\theta &=& Y_i^\theta -    f(t_i,X_i,\mu_i,Y_i^\theta) \Delta t_i  +  \Zc_{\tilde\theta}(t_i,\mu_i)(X_i) \Delta W_{t_i}, \quad i=0,\dots,N_T-1, 
\enqs
\end{footnotesize}
starting from $Y_0^\theta$ $=$ $\Uc_{\bar \theta}(\mu_0)(X_0)$,  from samples of $\mu_0$, and $X_0$ $\sim$ $\mu_0$.

\begin{Remark}
The algorithm in \cite{carlau19}  gives a solution to  the master equation only  for a given initial distribution while the algorithms presented here permit to solve the problem for all  initial distributions.
Furthermore, with  local algorithms, we are able to  obtain  the solution depending on $x$ and $\mu$ at each time step.
\end{Remark}

\vspace{2mm}

 \paragraph{Tests.}  For the numerical tests, we choose  $T=0.1$, $\kappa = 0.2$, $\sigma =0.5$,  $w(x)$ $=$ $\cos(x)$, and $a$ $=$ $0.1$  in \eqref{deff}. 
We take $\hat M=10$, $N=100000$, a number of bins of $K$ $=$ $200$ with a domain $\Kc$ $=$ $[-1.3,1.3]$ , $8E4$ gradient iterations at each optimization with an initial learning rate of $1E-3$ for the ADAM method. 
After optimization, we calculate the value function 
using the network at time $t=0$ (bin or cylindrical), and then estimate the associated  MSE for various test distributions $\mu^{test}$ as in  Section \ref{sec:numtest},  and following the local/global regression/BSDE  algorithms. 
The results \rd{ for the default hyperparameters of the networks (the same as for case A and B) } are reported in  Tables \ref{TableR} and  \ref{TableBSDE}, and show that the local BSDE algorithm with cylindrical neural networks provides  the best results.
\rd{ Notice in our case that the solution is very regular as it can be written as
$v(t,x,\mu)$ $=$ $e^{T-t}  [ \cos(x) \E_{\xi\sim\mu}[\cos(\xi)] + \sin(x) \E_{\xi\sim\mu}[\sin(\xi)]]$ and as shown in the first sections, a low number of layers, neurons is expected to be sufficient.
This is confirmed in table \ref{TableRRaf} and \ref{TableBSDERaf} where one hidden layer is added in the computations.
}
  
 \vspace{3mm}

 \begin{table}[H]
 \centering
 {
 \begin{tabular}{|c|c|c|c|c|}
 \hline
 Method & Network & Test 1 & Test 2 & Test 3 \\ 
 \hline
 Global & Bins &   1.55E-02  & 2.44E-02 & 1.46E-02   \\   
 Global & Cylinder & 6.70E-03  & 6.81E-03 & 6.04E-03  \\
 Local & Bins &  3.04E-02  & 2.32E-02 & 2.33E-02  \\   
 Local & Cylinder  & 1.11E-04  & 1.18E-04 & 2.35E-04    \\ 
 \hline
 \end{tabular}
 }
 \centering
 \caption{Regression  approach : MSE at time $0$ for PDE resolution (2 layers for the  cylinder network, 3 layers for the bin network).}  \label{TableR} 
 \end{table}

\begin{table}[H]
 \centering
 {
 \begin{tabular}{|c|c|c|c|c|} 
 \hline
 Method & Network &   Test 1 & Test 2 & Test 3  \\ 
 \hline
 Global & Bins &  2.53E-03  & 5.99E-03 & 1.96E-03    \\
 Global & Cylinder &  4.59E-04  & 3.81E-04 & 4.13E-04  \\
 Local & Bins &   2.92E-02  & 3.73E-02 & 2.36E-02 \\
 Local & Cylinder & 5.41E-05  & 6.45E-05 & 9.05E-05 \\
 \hline
 \end{tabular}
 }
 \centering
 \caption{BSDE approach : MSE at time $0$ for PDE resolution (2 layers for the  cylinder network, 3 layers for the bin network).} \label{TableBSDE} 
 \end{table}

\begin{table}[H]
 \centering
 {
 \begin{tabular}{|c|c|c|c|c|}
 \hline
 Method & Network & Test 1 & Test 2 & Test 3 \\ 
 \hline
 Global & Bins &  5.13E-02  & 8.15E-03 & 5.01E-02   \\   
 Global & Cylinder &  8.52E-04  & 1.07E-03 & 8.53E-04  \\
 Local & Bins &   4.09E-02  & 6.57E-02 & 3.16E-02  \\   
 Local & Cylinder  &  2.52E-03  & 1.96E-03 & 1.84E-03   \\ 
 \hline
 \end{tabular}
 }
 \centering
 \caption{Regression  approach : MSE at time $0$ for PDE resolution (3 layers for the  cylinder network, 4 layers for the bin network).}  \label{TableRRaf} 
 \end{table}

\begin{table}[H]
 \centering
 {
 \begin{tabular}{|c|c|c|c|c|} 
 \hline
 Method & Network &   Test 1 & Test 2 & Test 3  \\ 
 \hline
 Global & Bins &  1.74E-02  & 4.94E-03 & 1.39E-02   \\
 Global & Cylinder &  4.88E-03  & 3.66E-03 & 3.67E-03 \\
 Local & Bins &   3.92E-02  & 6.26E-02 & 3.14E-02  \\
 Local & Cylinder &  2.08E-03  & 1.70E-03 & 1.62E-03\\
 \hline
 \end{tabular}
 }
 \centering
 \caption{BSDE approach : MSE at time $0$ for PDE resolution (3 layers for the  cylinder network, 4 layers for the bin network).} \label{TableBSDERaf} 
 \end{table}

\begin{figure}[H]
     \centering
   \begin{minipage}[t]{0.32\linewidth}
  \centering
 \includegraphics[width=\textwidth]{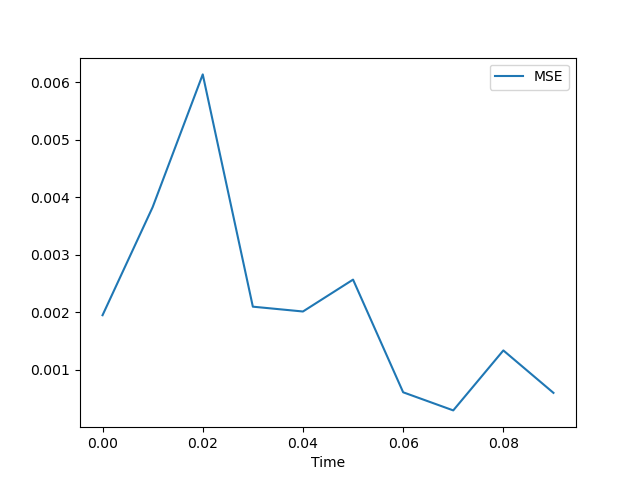}
\caption*{\tiny Test 1}
\end{minipage}
     \begin{minipage}[t]{0.32\linewidth}
  \centering
 \includegraphics[width=\textwidth]{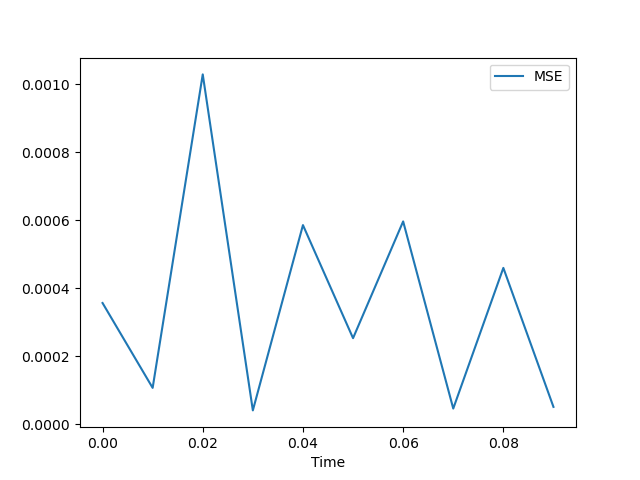}
\caption*{\tiny Test 2}
\end{minipage}
   \begin{minipage}[t]{0.32\linewidth}
  \centering
 \includegraphics[width=\textwidth]{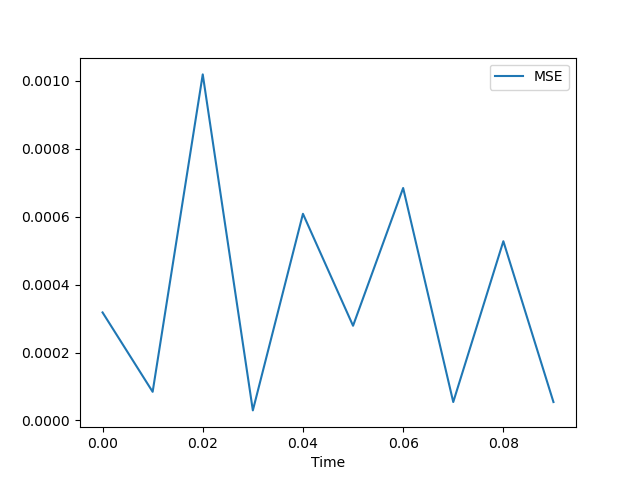}
\caption*{\tiny Test 3}
\end{minipage}
\caption{Local BSDE approach with cylindrical NN: MSE at different time steps.} \label{fig:localBSDEcyl} 
 \end{figure}

\vspace{3mm}
We plot in Figure \ref{fig:localBSDEcyl} the MSE error at different time steps when using the local BSDE algorithm with cylindrical mean-field neural networks. 
 
\begin{figure}[H]
     \centering
   \begin{minipage}[t]{0.32\linewidth}
  \centering
 \includegraphics[width=\textwidth]{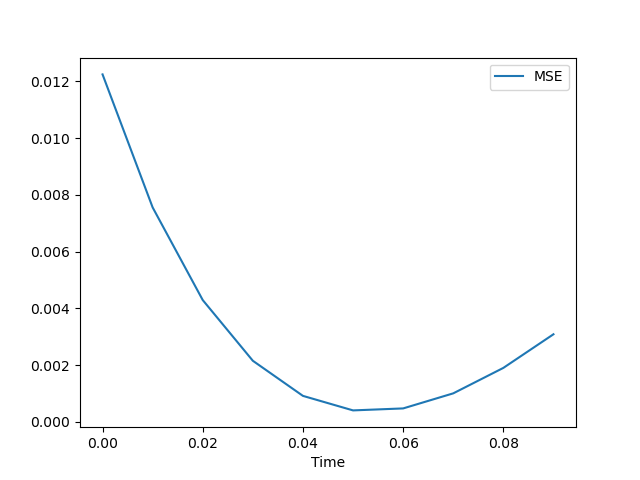}
\caption*{\small Test 1}
\end{minipage}
     \begin{minipage}[t]{0.32\linewidth}
  \centering
 \includegraphics[width=\textwidth]{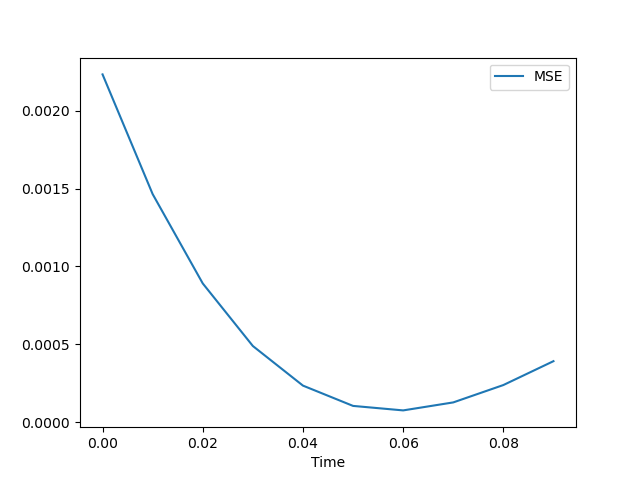}
\caption*{\small Test 2}
\end{minipage}
   \begin{minipage}[t]{0.32\linewidth}
  \centering
 \includegraphics[width=\textwidth]{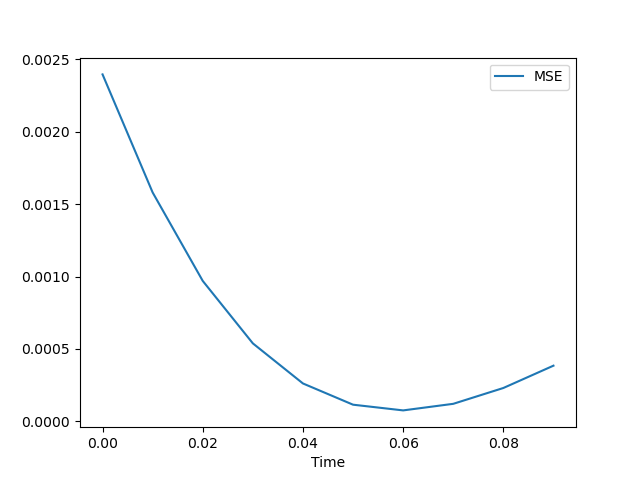}
\caption*{\small Test 3}
\end{minipage}
\caption{Global BSDE approach with cylindrical NN: MSE at different time steps.} \label{fig:globalBSDEcyl} 
 \end{figure}

Finally, we plot in Figure \ref{fig:globalBSDEcyl} the MSE error at different time steps when using the global BSDE algorithm with cylindrical mean-field neural networks. 

 \begin{Remark}
 The use of a single network in global algorithms permits  to have a smooth time representation of the value function as shown on Figure \ref{fig:globalBSDEcyl}.
 \end{Remark}

\section{Acknowledgements}
We thank Maximilien Germain and Mathieu Laurière for helpful discussions.
This work is supported by  FiME, Laboratoire de Finance des Marchés de l'Energie, and the ''Finance and Sustainable Development'' EDF - CACIB Chair.

\appendix

\section{Proofs of universal approximation theorems for mean-field neural networks} \label{sec:appen}

\subsection{Proof of Theorem \ref{theounivbin}}

Let $\varepsilon$ $>$ $0$ be given arbitrarily. Fix $\Kc$  a  bounded rectangular domain in  $\R^d$, and divide it into $K$ $=$ $\bar K^d$ bins:  Bin$(k)$, $k$ $=$ $1,\ldots,K$, of center $x_k$, and same area size 
$h$ $=$ $\lambda_d(\Kc)/K$, where $\lambda_d$ is the Lebesgue measure on $\R^d$.  Denote by  diam$_k$  the diameter of Bin$(k)$,  and notice that diam$_k$  $\leq$ ${\rm diam}(\Kc)/\bar K$, where 
diam$(\Kc)$ is the diameter of $\Kc$.   Let $V$ be  a continuous function on $\R^d\times\Pc_2(\R^d)$.

\vspace{1mm}

\noindent {\it Step 1.} 
For $\mu$ $\in$ $\Dc_c(\Kc)$  with density $\mrp^\mu$,  denote by $\hat\mu^K$ $=$ $\Lc_D(\bop^\mu)$  the probability measure with bin density $\bop^\mu$ in $\Dc_K$.  
Since $\mu$, $\hat\mu^K$ are supported on  the compact set $\Kc$, they lie in $\Pc_1(\R^d)$ the set of probability measures with finite first moment.
From the Kantorovich-Rubinstein dual representation of the $1$-Wasserstein distance, we have 
\beqs
\Wc_1(\mu,\hat\mu^K) & \leq & \sup_{\phi} \int_{\Kc}  \phi(x) (\mrp^\mu(x) - \hat\mrp_\Kc^\mu(x)) \d x,  
\enqs
where $\hat\mrp_\Kc^\mu(x)$ $=$ $\sum_{k=1}^K \frac{\mrp^\mu(x_k)}{N_K} 1_{x \in {\rm Bin}(k)}$, with $N_K$ $=$ $\sum_{k=1}^K \mrp^\mu(x_k) h$,  and the supremum is taken over all Lipschitz continuous functions $\phi$ on $\Kc$ with Lipschitz constant bounded by $1$, and where we can assume w.l.o.g. that $\phi(x_0)$ $=$ $0$ for some fixed point $x_0$ in $\Kc$.  We then have 
\beqs
\Wc_1(\mu,\hat\mu^K) & \leq &  \sup_{\phi} \sum_{k=1}^K \int_{{\rm Bin}(k)} \phi(x) \big( \mrp^\mu(x) - \frac{\mrp^\mu(x_k)}{N_K} \big) \d x \\
& \leq &   \sup_{\phi} \sum_{k=1}^K \Big[ \int_{{\rm Bin}(k)} \phi(x) \big( \mrp^\mu(x) - \mrp^\mu(x_k) \big) \d x \\
& & +  \int_{{\rm Bin}(k)} \phi(x)  \mrp^\mu(x_k) \frac{N_K-1}{N_K} \d x \Big] \\
& \leq & {\rm diam}(\Kc)  h  \sum_{k=1}^K \omega_\Kc^\mu( {\rm diam}_k)    +  {\rm diam}(\Kc) |N_K - 1| \\
& \leq & 2  {\rm diam}(\Kc) \lambda_d(\Kc)  \omega_\Kc^\mu \big(  \frac{{\rm diam}(\Kc)}{K^{\frac{1}{d}}}\big), 
\enqs
where we used in the third inequality the fact that $|\phi(x)|$ $\leq$ $|x-x_0|$ $\leq$ ${\rm diam}(\Kc)$, for any $x$ $\in$ $\Kc$, and $|\mrp^\mu(x)-\mrp^\mu(x_k)|$ $\leq$  $\omega_\Kc^\mu({\rm diam}_k)$ for any $x$ $\in$ ${\rm Bin}(k)$, and in the fourth inequality the fact  that diam$_k$  $\leq$ ${\rm diam}(\Kc)/K^{\frac{1}{d}}$,  $Kh$ $=$ $\lambda_d(\Kc)$, and the relation 
\beqs
\big| 1  -N_k \big| &=& 
\big| \sum_{k=1}^K \int_{{\rm Bin}(k)} \big[ \mrp^\mu(x) - \mrp^\mu(x_k) \big] \d x \Big| 
\; \leq \;   \lambda_d(\Kc)  \omega_\Kc^\mu\big(\frac{ {\rm diam}(\Kc) }{K^{\frac{1}{d}}}\big). 
\enqs
By noting that $\Wc_2(\mu,\hat\mu^K)$ $\leq$ $\sqrt{{\rm diam}(\Kc)\Wc_1(\mu,\hat\mu^K)}$,
this shows that 
\begin{align} \label{convWK} 
\sup_{\mu \in \Dc_c^{\bar\omega}(\Kc)} \Wc_2(\mu,\hat\mu^K)  & \rightarrow \; 0, \quad \mbox{ as } K \rightarrow \infty. 
\end{align} 
On the other hand, by Lemma 5.7 and Proposition 5.3 in \cite{car13}, the set $\Dc_c(\Kc)$ is relatively compact in $\Pc_2(\R^d)$, and thus $V$ is uniformly continuous on $\Kc\times\Dc_c(\Kc)$. From \eqref{convWK}, it follows that  there exists 
$K$ $\in$ $\N^*$, such that 
\begin{align} \label{inegV1} 
| V(x,\mu) - V(x,\hat\mu^K) | & \leq \; \frac{\varepsilon}{2}, \quad \forall  x \in \Kc, \; \mu \in \Dc_c^{\bar\omega}(\Kc). 
\end{align}

\vspace{1mm}

\noindent {\it Step 2.} 
Denote by $V_K$ the function defined on $\R^d\times\Dc_K$ by 
\beqs
V_K(x,\bop) &=& V(x,\Lc_D(\bop)),  \quad (x,\bop) \in \R^d\times\Dc_K,
\enqs
where we recall that $\Dc_K$ $=$ $\{ \bop = (p_k)_{k\in\llbracket 1,K\rrbracket} \in \R_+^K: \sum_{k=1}^K p_k h =1\}$, and $\Lc_D(\bop)$ is the probability measure with bin density $\bop$ $\in$ $\Dc_K$.  It is clear that when $(\bop^n)_n$ converges to $\bop$ in $\Dc_K$, 
then $\Lc_D(\bop^n)$ converges weakly to $\Lc_D(\bop)$, and thus $V_K$ is continuous on $\R^d\times\Dc_K$.  By the classical universal approximation theorem for finite-dimensional functions (see \cite{hor91}), there exists a feedforward neural network $\Phi$ on 
$\R^d\times\Dc_K$ such that 
\begin{align} \label{inegV2} 
| V_K(x,\bop) - \Phi(x,\bop) | & \leq \; \frac{\varepsilon}{2}, \quad \forall x \in \Kc, \; \bop \in \Dc_K. 
\end{align}
We conclude that for all $x$ $\in$ $\Kc$, $\mu$ $\in$ $\Dc_c^{\bar\omega}(\Kc)$, 
\begin{align}
\big| V(x,\mu) -  \Phi(x,\bop^\mu) \big| & \leq \;  \big| V(x,\mu) -  V(x,\hat\mu^K) \big| + \big| V_K(x,\bop^\mu) -  \Phi(x,\bop^\mu) \Big| \; \leq \; \varepsilon.
\end{align}
by noting that $V(.,\hat\mu^K)$ $=$ $V_K(.,\bop^\mu)$, and using \eqref{inegV1}-\eqref{inegV2}.

\subsection{Proof of Theorem \ref{theounivcyl2}}

{\it Step 1: Approximation theorem on compact set}.  Let $\varepsilon$ $>$ $0$ be given arbitrarily. 
Fix $\Kc$ a compact set of $\R^d$, and let $V$ be  a continuous function on $\R^d\times\Pc_2(\R^d)$.  By the density of cylindrical polynomial function with respect to mean-field functions, see Lemma 3.12 in \cite{guophawei22}, for all $\varepsilon$ $>$ $0$, there exists $k$ 
$\in$ $\N^*$,  $P$ a polynomial function from $\R^d\times\R^k$ into $\R$, $Q$ a polynomial function from $\R^d$ into $\R^k$ s.t.  
\begin{align}
\big| V(x,\mu) - P(x,<Q,\mu>) \big| & \leq \;  \frac{\varepsilon}{3}, \quad \forall x \in \Kc, \; \mu \in \Pc(\Kc). 
\end{align}
Now, by the uniform continuity of $P$ on compact sets of $\R^d\times\R^k$, and the classical universal approximation theorem for finite-dimensional functions applied to $Q$, there exists a feedforward neural network $\varphi$ from $\R^d$ into $\R^k$ such that 
\begin{align}
\big| P(x,<Q,\mu>) - P(x,<\varphi,\mu>) \big| & \leq \;  \frac{\varepsilon}{3}, \quad \forall x \in \Kc, \mu \in \Pc(\Kc), 
\end{align}
by noting that one can find some compact set $\Yc$ (depending on $Q$ and $\Kc$) of $\R^k$ such that  $<Q,\mu>$ and then $<\varphi,\mu>$  lie in  $\Yc$ for all $\mu$ $\in$ $\Kc$. Next, we invoke again  the  classical universal approximation theorem for finite-dimensional functions to get the existence of  a feedforward neural network $\Psi$ on $\R^d\times\R^k$ such that 
\begin{align}
\big| P(x,y) - \Psi(x,y) \big| & \leq \;  \frac{\varepsilon}{3}, \quad \forall  (x,y) \in \Kc\times\Yc. 
\end{align}
We conclude that for all $x$ $\in$ $\Kc$, $\mu$ $\in$ $\Pc(\Kc)$, 
\begin{align} 
& \;  \big| V(x,\mu) - \Psi(x,<\varphi,\mu>) \big| \nonumber \\
& \leq \;     \big| V(x,\mu) - P(x,<Q,\mu>) \big|  +     \big| P(x,<Q,\mu>) - P(x,<\varphi,\mu>) \big| \nonumber \\
&  \quad  \; +   \;  \big| P(x,<\varphi,\mu>) - \Psi(x,<\varphi,\mu>) \big|  \; \leq \; \varepsilon.  \label{Vpsi} 
\end{align} 

 \vspace{2mm}


\noindent  {\it Step 2:  Approximation theorem in $L^2$}.  Let $\varepsilon$ $>$ $0$ be given arbitrarily, and $\nu$ be a probability measure on $\Pc_2(\R^d)$. 
 Given $M$ $>$ $0$, we truncate the function $V$ by defining $V_M$ on $\R^d\times\Pc_2(\R^d)$ as
 \begin{equation*}
 V_M(x,\mu) \;  =  \; \left\{ 
 \begin{array}{cc} 
 V(x,\mu), & \quad  \mbox{ if } \;  |V(x,\mu)|  \leq M,  \\
 M \frac{V(x,\mu)}{|V(x,\mu)|}, & \quad  \mbox{ if } \; |V(x,\mu)| > M,
 \end{array}
 \right. 
 \end{equation*}
 so that $ |V_M(x,\mu)| \leq M$ for all $x$ $\in$ $\R^d$, $\mu$ $\in$ $\Pc_2(\R^d)$.  It is clear that $V_M$ converges pointwise to $V$ as $M$ goes to infinity, and thus by the dominated convergence theorem $\|V -V_M\|^2_{L^2(\nu)}$ converges to zero. We can 
 thus choose $M$ $>$ $\frac{\sqrt{\varepsilon}}{4}$ so that 
 \begin{align} \label{VMV} 
 \|V -V_M\|^2_{L^2(\nu)} & = \; \int_{\Pc_2(\R^d)} |V(.,\mu) - V_M(.,\mu)|_\mu^2 \nu(\d \mu) \; \leq \; \frac{\varepsilon}{8}. 
 \end{align} 
 
 Next, we consider some compact set $\Kc$ of $\R^d$ such that $\nu(\Pc_2(\R^d)\setminus\Pc(\Kc))$ $\leq$ $\varepsilon/(80 M^2)$, and we note that $V_M$ is continuous  on $\R^d\times\Pc_2(\R^d)$. We can then apply the universal approximation theorem in 
 Step 1,  to get the existence of 
 a cylindrical mean-field neural network $(\tilde\Psi: \R^d\times\R^k \mapsto \R^p,\varphi:\R^d\mapsto\R^k)$ s.t. 
\begin{align} \label{VMcyl1} 
\big| V_M(x,\mu) -  \tilde\Psi(x,<\varphi,\mu>) \big| & \leq \; \frac{\sqrt{\varepsilon}}{4}, \quad \forall x \in \Kc, \; \mu \in \Pc(\Kc).  
\end{align}
This implies in particular that
\begin{align}
\big|\tilde\Psi(x,<\varphi,\mu>) \big|  & \leq \;   | V_M(x,\mu)| +   \big| V_M(x,\mu) - \tilde\Psi(x,<\varphi,\mu>) \big| \nonumber \\
&  \leq \;   M + \frac{\sqrt{\varepsilon}}{4}  < 2M, \quad \quad x \in \Kc, \; \forall \mu \in \Pc(\Kc).  \label{Psibor}
\end{align} 
By the clipping lemma C.1 in \cite{lanetal21}, there exists a neural network $\gamma$ $:$ $\R^p$ $\mapsto$ $\R^p$, such that 
\begin{equation} \label{clip} 
\left\{
\begin{array}{rcl} 
|\gamma(y) - y | & \leq  & \frac{\sqrt{\varepsilon}}{4}, \quad \mbox{ if } |y| \leq  M + \frac{\sqrt{\varepsilon}}{4}, \\
|\gamma(y)| & \leq &  2M, \quad \forall y \in \R^p. 
\end{array}
\right.
\end{equation}
(Actually, when $p$ $=$ $1$, one can simply take $\gamma(y)$ $=$ $\min[\max[y,-2M],2M]$).  Define now the neural network $\Psi$ $:$ $\R^d\times\R^k$ $\rightarrow$ $\R^p$ by 
$\Psi$ $=$ $\gamma\circ\tilde\Psi$.  It is then bounded from  above by 
\begin{align} \label{Psibor2M} 
|\Psi(x,z)| & \leq \; 2M, \quad \forall (x,z) \in \R^d\times\R^k,
\end{align} 
and satisfies  for all $x$ $\in$ $\Kc$, $\mu$ $\in$ $\Pc(\Kc)$, 
\begin{align}
\big| V_M(x,\mu) - \Psi(x,<\varphi,\mu>) \big| & \leq \; \big| V_M(x,\mu) - \tilde \Psi(x,<\varphi,\mu>) \big| \\
& \quad + \; \big| \gamma \circ \tilde\Psi(x,<\varphi,\mu>)  -   \tilde \Psi(x,<\varphi,\mu>) \big| \\
& \leq \; \frac{\sqrt{\varepsilon}}{4} +  \frac{\sqrt{\varepsilon}}{4} \; = \;  \frac{\sqrt{\varepsilon}}{2},
\end{align}
by \eqref{VMcyl1}, \eqref{Psibor}, and \eqref{clip}.  It follows  that 
\begin{align}
 C & :=  \int_{\Pc_2(\R^d)} | V_M(.,\mu) -  \Psi(.,<\varphi,\mu>) \big|_\mu^2  \nu(\d \mu)  \\  & \leq    \;  \int_{\Pc(\Kc)} | V_M(.,\mu) -  \Psi(.,<\varphi,\mu>) \big|_\mu^2  \nu(\d \mu)  \\
& \; +  \;  2 \int_{\Pc_2(\R^d)\setminus\Pc(\Kc)} ( |V_M(.,\mu)|_\mu^2  + | \Psi(.,<\varphi,\mu>)|_\mu^2 ) \nu(\d \mu) \\
& \leq \;  \frac{\varepsilon}{4} + 2(M^2 + 4M^2) \frac{\varepsilon}{80 M^2}  \; = \; \frac{3\varepsilon}{8}. 
\end{align} 

We conclude with \eqref{VMV} that  
\begin{align}
D := & \int_{\Pc_2(\R^d)} | V(.,\mu) -  \Psi(.,<\varphi,\mu>) \big|_\mu^2  \nu(\d \mu) \\
& \leq \:  2 \int_{\Pc_2(\R^d)} |V(.,\mu) - V_M(.,\mu)|_\mu^2 \nu(\d \mu) \\
& \quad + \;  2 \int_{\Pc_2(\R^d)} | V_M(.,\mu) -  \Psi(.,<\varphi,\mu>) \big|_\mu^2  \nu(\d \mu)   \\
& \leq \; 2 \frac{\varepsilon}{8} +  2 \frac{3\varepsilon}{8}  \; = \; \varepsilon. 
\end{align}

\bibliographystyle{elsarticle-num} 
\bibliography{bibli}

\end{document}